\documentclass[12pt]{amsart}
\usepackage{amsmath}
\usepackage{amsxtra}
\usepackage{amscd}
\usepackage{amsthm}
\usepackage{amsfonts}
\usepackage{amssymb}
\usepackage{eucal}

\input prepictex
\input pictex
\input postpictex

\newtheorem{thm}{Theorem}[section]
\newtheorem{lem}{Lemma}[section]
\newtheorem{cor}{Corollary}[section]
\newtheorem{prop}{Proposition}[section]

\theoremstyle{definition}

\theoremstyle{remark}
\newtheorem{rem}{Remark}[section]

\numberwithin{equation}{section}

\begin{document}

\newcommand{\thmref}[1]{Theorem~\ref{#1}}
\newcommand{\secref}[1]{Section~\ref{#1}}
\newcommand{\lemref}[1]{Lemma~\ref{#1}}
\newcommand{\propref}[1]{Proposition~\ref{#1}}
\newcommand{\corref}[1]{Corollary~\ref{#1}}
\newcommand{\remref}[1]{Remark~\ref{#1}}
\newcommand{\eqnref}[1]{(\ref{#1})}
\newcommand{\exref}[1]{Example~\ref{#1}}

\newcommand{\nc}{\newcommand}
\nc{\on}{\operatorname} \nc{\Z}{{\mathbb Z}} \nc{\C}{{\mathbb C}}
\nc{\oo}{{\mf O}}

\nc{\bib}{\bibitem}
\nc{\pa}{\partial}
\nc{\F}{{\mf F}}
\nc{\rarr}{\rightarrow}
\nc{\larr}{\longrightarrow}
\nc{\al}{\alpha}
\nc{\ri}{\rangle}
\nc{\lef}{\langle}
\nc{\W}{{\mc W}}
\nc{\gam}{\ol{\gamma}}
\nc{\Q}{\ol{Q}}
\nc{\q}{\widetilde{Q}}
\nc{\la}{\lambda}
\nc{\ep}{\epsilon}
\nc{\g}{\mf g}
\nc{\h}{\mf h}
\nc{\n}{\mf n}
\nc{\A}{{\mf a}}
\nc{\G}{{\mf g}}
\nc{\Li}{{\mc L}}
\nc{\La}{\Lambda}
\nc{\is}{{\mathbf i}}
\nc{\V}{\mf V}
\nc{\bi}{\bibitem}
\nc{\NS}{\mf N}
\nc{\dt}{\mathord{\hbox{${\frac{d}{d t}}$}}} \nc{\E}{\mc E}
\nc{\ba}{\tilde{\pa}}
\def\smapdown#1{\big\downarrow\rlap{$\vcenter{\hbox{$\scriptstyle#1$}}$}}
\nc{\mc}{\mathcal} \nc{\mf}{\mathfrak} \nc{\ol}{\fracline}
\nc{\el}{\ell}
\nc{\etabf}{{\bf \eta}}
\nc{\x}{{\bf x}}
\nc{\xibf}{{\bf \xi}}
\nc{\y}{{\bf y}}
\advance\headheight by 2pt

\title{ Howe Duality for Lie Superalgebras}

\author[Shun-Jen Cheng]{Shun-Jen Cheng$^*$}
\thanks{$^*$partially supported by NSC-grant 89-2115-M-006-002 of the R.O.C}
\address{Department of Mathematics, National Taiwan University, Taipei,
Taiwan} \email{chengsj@math.ntu.edu.tw}

\author{Weiqiang Wang}
\address{Department of Mathematics, North Carolina State University,
Raleigh, NC 27695-8205, USA}
\email{wqwang@math.ncsu.edu}

\begin{abstract}
We study a dual pair of general linear Lie superalgebras in the
sense of R.~Howe. We give an explicit multiplicity-free
decomposition of a symmetric and skew-symmetric algebra (in the
super sense) under the action of the dual pair and present
explicit formulas for the highest weight vectors in each isotypic
subspace of the symmetric algebra. We give an explicit
multiplicity-free decomposition into irreducible $gl(m|n)$-modules
of the symmetric and skew-symmetric algebras of the
symmetric square of the natural representation of $gl(m|n)$. In
the former case we find as well explicit formulas for the highest
weight vectors. Our work unifies and generalizes the classical
results in symmetric and skew-symmetric models and admits several
applications.
\vspace{.3cm}

\noindent{\bf Key words:} Lie superalgebra, Howe duality,
highest weight vectors.

\vspace{.3cm}

\noindent{\bf Mathematics Subject Classifications (1991)}: 17B67.

\end{abstract}

\maketitle

\section{Introduction}

Howe duality is a way of relating representation theory of a pair
of reductive Lie groups/algebras \cite{H1, H2}. It has found many
applications to invariant theory, real and complex reductive
groups, $p$-adic groups and infinite-dimensional Lie algebras etc.

As an example we consider one of the fundamental cases---the
$(gl(m), gl(n) )$ Howe duality. The symmetric algebra $S
(\C^m\otimes\C^n)$ and the skew-symmetric algebra $\Lambda
(\C^m\otimes\C^n)$ admit remarkable multiplicity-free
decompositions under the natural actions of $gl(m) \times gl(n)$.
The highest weight vectors of $gl(m) \times gl(n)$ inside the
symmetric algebra are given by products of certain determinants
(see (\ref{deltar})) and form a free abelian semi-group while
those inside the skew-symmetric algebra are given by Grassmann
monomials, cf. \cite{H2, KV, GW}.

Our present paper is devoted to the study of Howe duality for Lie
superalgebras and its applications. It is by now a well
established fact that one should put the Grassmann variables on
the same footing as Cartesian variables and hence it is natural to
consider the supersymmetric algebra, which is a mixed tensor of
symmetric and skew-symmetric algebras. In this paper we give a
complete description of the Howe duality in a
symmetric\footnote{In this paper we will freely suppress the term
{\em super}. So in case when a superspace is involved, the terms
symmetric, commute etc.~mean {\em super}symmetric, {\em
super}commute etc unless otherwise specified.} algebra under the
action of a dual pair of general linear Lie superalgebras and find
explicit formulas for the highest weight vectors inside our
symmetric model. A dual pair consisting of a general linear Lie
superalgebra and a general linear Lie algebra was discussed in
\cite{H1}. We also study in detail some other multiplicity-free
actions of the general linear Lie superalgebras as specified
below.

Our motivation is manifold. First, our work is motivated by an
attempt to unify the Howe duality in the symmetric and
skew-symmetric models \cite{H2} which have many differences and
similarities. Specialization of our results gives rise to the Howe
duality for general linear Lie algebras in both symmetric and
skew-symmetric models. Secondly, we are motivated by our study of
the duality in the infinite-dimensional setup (see the review
\cite{W} and references therein) and our work in progress on its
generalization to the superalgebra case. We realize that we have
to understand the finite-dimensional picture better first in order
to have a more complete description of the infinite-dimensional
picture. Thirdly, there exists a new type of Howe duality which is
of pure superalgebra phenomenon which is treated in \cite{CW}.

Let us discuss the contents of the paper. The generalization of
Schur duality for the superspace was given by Sergeev
\cite{Se}. For lack of an analog in the super setup of the
criterion of multiplicity-free action in terms of the existence of
a dense open orbit of a Borel subgroup (cf. \cite{V} and
\cite{H2}), we use Sergeev's result to derive the decomposition of
the symmetric algebra $S(\C^{p|q}\otimes\C^{m|n})$ with respect to
the action of the sum of two general linear Lie superalgebras
$gl(p|q)\times gl(m|n)$. We see that a representation of $gl(p|q)$
is paired with a representation of $gl(m|n)$ parameterized by the
same Young diagram. On the other hand one can show that our Howe
duality for superalgebras implies Sergeev's Schur duality as well.
We also obtain an explicit multiplicity-free decomposition of a
skew-symmetric algebra $\Lambda(\C^{p|q}\otimes\C^{m|n})$ as
$gl(p|q)\times gl(m|n)$-modules. In particular it follows that
$gl(p|q)$ and $gl(m|n)$ when acting on
$S(\C^{p|q}\otimes\C^{m|n})$ and respectively on
$\Lambda(\C^{p|q}\otimes\C^{m|n})$ are mutual (super)centralizers.
A remarkable phenomenon is the complete reducibility of the
symmetric model under the action of the dual
pair, which is quite unusual for Lie superalgebras.

In a purely combinatorial way, Brini, Palareti and Teolis
\cite{BPT} were indeed the first to obtain an explicit
decomposition of $S(\C^{p|q}\otimes\C^{m|n})$ under the action of
$gl(p|q)\times gl(m|n)$. In addition, their combinatorial approach
exhibits explicit bases parameterized by so-called left (or right)
symmetrized bitableuax between two ``standard Young diagrams''
(see \cite{BPT} for definition). However, Brini {\em et al} did
not identify the highest weights for these $gl(p|q)\times
gl(m|n)$-modules inside $S(\C^{p|q}\otimes\C^{m|n})$.

We also obtain an explicit decomposition into irreducible
$gl(m|n)$-modules of the symmetric algebra $S({S^2\C^{m|n}})$ and
respectively skew-symmetric algebra $\Lambda({S^2\C^{m|n}})$ of
the symmetric square of the natural representation of $gl(m|n)$.
These results unifies and generalizes several classical results
and they can be proved in an analogous way as in the classical
case \cite{H2, GW}.

Associated to the Howe duality and the above $gl(m|n)$-module
decompositions, we obtain, by taking characters, various
combinatorial identities involving the so-called hook Schur
functions. Being generalization of Schur functions, these hook
Schur functions have been studied in \cite{BR}. The
decompositions mentioned above in turn provide the representation
theoretic realization of the corresponding combinatorial
identities. For example, the $(gl(m|n), gl(p|q))$-duality gives
rise to a combinatorial identity for the hook Schur functions
which generalizes the Cauchy identity for Schur functions, cf.
\cite{H2}. Specializations and variations of these combinatorial
identities are well known and other proofs can be found in
Macdonald \cite{M}.

However it is a much more difficult problem to find explicit
formulas for the highest weight vectors of $gl(p|q)\times
gl(m|n)$-modules inside the symmetric algebra. We first find
formulas for the highest weight vectors in the case for $q =0$
(and so for $n=0$ by symmetry). A main ingredient in the formulas
for the highest weight vectors is given by the determinant of a
matrix which involves both Cartesian variables $x_j^i$'s and
Grassmann variables $\eta^i$'s of the form:
\begin{equation*}
\begin{pmatrix}
x_1^1&x_1^2&\cdots &x_1^r\\ x_2^1&x_2^2&\cdots &x_2^r\\
\vdots&\vdots&\cdots &\vdots\\ x_m^1&x_m^2&\cdots &x_m^r\\
 \eta^1&\eta^2&\cdots &\eta^r\\
\vdots&\vdots&\cdots &\vdots\\ \eta^1&\eta^2&\cdots &\eta^r\\
\end{pmatrix}.
\end{equation*}
We remark that the rows involving Grassmann variables are the same
but the determinant is nonzero (one needs to overcome some
psychological barriers). Note that when $m =r$ the Grassmann
variables disappear and the above determinant reduces to those
mentioned earlier which occur in the formulas for highest weight
vectors in the symmetric algebra case of the classical Howe
duality. When $m =0$, the Cartesian variables disappear and the
above determinant is equal to (up to a scalar multiple) a
Grassmann monomial which shows up in the formulas for highest
weight vectors in the classical skew-symmetric algebra case.

We show that the $gl(p|q)\times gl(m|n)$ highest weight vectors
form an abelian  semigroup in the case when $p =m$. However in
contrast to the Lie algebra case this semigroup is not free in
general. We find that the generators of the semigroup are given by
highest weight vectors associated to rectangular Young diagrams of
length not exceeding $m +1$. This way we are able to find explicit
formulas for all highest weight vectors in the case when $q =0$
(or $m =0$), or $p =m$.

In the general case the highest weight vectors no longer form a
semigroup. We find a nice way to overcome this difficulty by
introducing some extra variables which, roughly speaking, help us
to reduce the general case to the case $p=m$. Then we use a simple
method to get rid of the extra variables to obtain the genuine
highest weight vectors we are looking for.

In contrast to the Howe duality in the Lie algebra setup, it is
difficult to check directly the highest weight condition of the
vectors we have obtained. We use instead the multiplicity-free
decomposition of the symmetric algebra to get around this
difficulty. As highest weight constraint we obtain interesting
non-trivial polynomial identities typically involving various
minors of a matrix.

We also find explicit formulas for the highest weight vectors
appearing in the $gl(m|n)$-module decomposition of
$S({S^2\C^{m|n}})$. These highest weight vector formulas, which
constitute a mixture of determinants and Pfaffians, have somewhat
similar features as those found in the Howe duality for the
general linear Lie superalgebras.

A formula for highest weight vectors in the decomposition of
$S(\C^{p|q}\otimes\C^{m|n})$ as $gl(p|q)\times gl(m|n)$ modules
may also be obtained in principle using the combinatorial approach
of Brini {\em et al} \cite{BPT}, and in this way the highest
weights for these $gl(p|q)\times gl(m|n)$ modules can be
identified. However this way one can neither expect to obtain
formulas as explicit as ours, nor can one see the semigroup structure
of the set of the highest weight vectors which is the guiding
principle for us to find these vectors. It is also interesting to
see whether our results concerning the decomposition of
$S(S^2\C^{m|n})$ (and respectively $\Lambda(S^2\C^{m|n})$) and the
highest weight vectors in these models may also be obtained with
extra insights from the combinatorial approach in \cite{BPT} as
well.

The plan of the paper goes as follows. In Section
\ref{sect_classical} we review the classical dual pairs of general
linear Lie algebras and Schur duality. In Section
\ref{sect_superdual} we present various multiplicity-free actions
for Lie superalgebras and obtain the corresponding symmetric
function identities. \secref{jointhwv} is devoted to the
construction of the $gl(p|q)\times gl(m|n)$ highest weight vectors
inside $S(\C^{p|q}\otimes\C^{m|n})$. More precisely, in Section
\ref{sect_hwto}, Section \ref{sect_pmequal}, and Section
\ref{sect_general}, we find explicit formulas of highest weight
vectors in the case $q =0$, $p =m$, and the general case,
respectively. Finally in \secref{hwvector} we construct the
$gl(m|n)$ highest weight vectors inside $S(S^2\C^{m|n})$.

{\bf Acknowledgment.} Our work is greatly influenced by the
beautiful article \cite{H2} of R.~Howe to whom we are grateful.
The results in this paper and and its sequel \cite{CW} are based
on our two preprints under the same title (part one and two) as
the current paper. After we submitted part one and finished part
two, we came across two preprints of Sergeev: ``An analog of the
classical invariant theory for Lie superlagebras", I, II,
math.RT/9810113 and math.RT/9904079, which have overlaps with our
work. More precisely, Sergeev obtained independently the
$gl(p|q)\times gl(m|n)$-module decomposition of the symmetric
algebra $S(\C^{p|q}\otimes\C^{m|n})$ (i.e. our
\thmref{supergl-duality}), and the $gl(m|n)$-module decomposition
of $S({S^2\C^{m|n}})$ (i.e. our \thmref{symmetric}). Finally we
also like to thank the referee for bringing the paper \cite{BPT},
which is discussed in the introduction, to our attention.
Therefore we have made corresponding changes on our preprints
which result in the current version of this paper and \cite{CW}.
\section{The Classical Picture}  \label{sect_classical}

In this section we will review some classical multiplicity-free
actions of the general linear Lie algebra. We begin with the
classical $gl(m)\times gl(n)$-duality, cf. Howe \cite{H2}.

Let $\la =(\la_1, \la_2, \ldots, \la_l)$ be a partition of the
integer $|\la|=\la_1 +\ldots +\la_l$, where $\la_1\geq \dots \geq
\la_l > 0$. The integer $|\la|$ is called the {\em size}, $l$ is
called the {\em length} (denoted by $l (\la )$), and $\la_1$ is
called the {\em width} of the partition $\la $. Let $\la'$ denote
the Young diagram obtained from $\la$ by transposing. We will often
denote $\lambda_1$ by $t$ and write $\la ' =
(\la_1',\la_2',\ldots,\la_t')$. For example, the Young diagram
\begin{equation}\label{standard}
{\beginpicture
\setcoordinatesystem units <1.5pc,1.5pc> point at 0 2
\setplotarea
x from 0 to 1.5, y from 0 to 4
\plot 0 0 0 4 1 4 1 0 0 0 /
\plot 1 1 2 1 2 4 1 4 /
\plot 2 2 3 2 3 4 2 4 /
\plot 3 3 5 3 5 4 3 4 /
\plot 0 1 1 1 /
\plot 0 2 2 2 /
\plot 0 3 3 3 /
\plot 4 3 4 4 /
\put{$r\left\{\hbox to 0pt{ \vrule width 0 pc height 2.5pc depth
2.5pc\hfil}
  \right.$} at -.5 2
\put{$\underbrace{\hbox to 8pc{\ \hfill}}_{t}$} at 2.5 -0.5
\endpicture}
\end{equation}
stands for the partition $(5,3,2,1)$ and its transpose is
the partition $(4,3,2,1,1)$.

Given a partition $\la =(\la_1, \la_2, \ldots, \la_l)$ satisfying
$l \leq m$, we may regard $\la$ as a highest weight of $gl(m)$ by
identifying $\la$ with the $m$-tuple $(\la_1, \la_2, \ldots,
\la_l, 0, \ldots, 0)$ by adding $m -l$ zeros to $\la$. We denote the
irreducible finite-dimensional highest weight module of $gl(m)$ by
$V^{\la}_m$.

Consider the natural action of the complex general linear Lie groups
$GL(m)$ and $ GL(n)$ on the space $\C^m\otimes\C^n$. If we
identify $\C^m\otimes\C^n$ with $M_{mn}$, the space of all $m
\times n$ matrices, then the actions of $GL(m)$ and $GL(n)$ are
given by left and right multiplications:
\[
(g_1, g_2) (T) = (g_1^t)^{-1} T g_2^{-1} \quad g_1 \in GL(m), g_2
\in GL(n), T \in M_{mn}.
\]
The Lie algebras $gl(m)$ and $gl(n)$ act on $\C^m\otimes\C^n$
accordingly. Denoting by $S (\C^m\otimes\C^n)$ the symmetric
tensor algebra of $\C^m\otimes\C^n$ with an induced action of
$gl(m)\times gl(n)$, we have the following multiplicity-free
decomposition of $S(\C^m\otimes\C^n)$ as a $gl(m)\times
gl(n)$-module:

\begin{equation*}
S(\C^m\otimes\C^n) \cong \sum_{\lambda} V_m^{\lambda}\otimes
V_n^{\lambda},
\end{equation*}
where the sum above is over Young diagrams $\lambda$ of length not
exceeding ${\rm min}(m,n)$.

One can find an explicit formula for the $gl(m)\times gl(n)$
highest weight vectors in this decomposition. Let us denote a
basis of $\C^m$ by $x_1,x_2,\ldots,x_m$ and a basis of $\C^n$ by
$x^1,x^2,\ldots,x^n$.  Then the vectors $x^j_i:=x_i \otimes x^j$,
for $i =1,\ldots, m$ and $j =1,\ldots, n$ form a basis for
$\C^m\otimes\C^n$ so that we may identify $S(\C^m\otimes\C^n)$
with $\C[x_{1}^1,\ldots,x_1^n,\ldots,x_m^1,\ldots,x_m^n]$. Using
this identification the standard Borel subalgebra of $gl(m)$ is a
sum of the Cartan subalgebra generated by
\[
\sum_{j=1}^n x_{i}^j\frac{\partial}{\partial x_i^j}, \quad 1\le
i\le m,
\]
and the nilpotent radical generated by
\[
\sum_{j=1}^n x_{i-1}^{j}\frac{\partial}{\partial x_{i}^j}, \quad
2\le i\le m.
\]
Similarly the Borel subalgebra of $gl(n)$ is the sum of the Cartan
subalgebra generated by
\[
\sum_{i =1}^m x_{i}^j\frac{\partial}{\partial x_i^j}, \quad 1\le j
\le n,
\]
and the nilpotent radical generated by
\[
\sum_{i =1}^m x_i^{j -1}\frac{\partial}{\partial x_i^j}, \quad
2\le j \le n.
\]

Let $\epsilon_i$ for $i=1,\ldots,m$ (respectively $\tilde{\epsilon}_j$
for $j=1,\ldots,n$) be the fundamental weights corresponding to
the Cartan subalgebra of $gl(m)$ (respectively $gl(n)$) above. For $1\le
r\le {\rm min}(m,n)$ define
\begin{equation}\label{deltar}
\Delta_r:={\rm det}\begin{pmatrix}
x_{1}^{1}&x_{2}^{1}&\cdots&x_{r}^{1}\\
x_{1}^{2}&x_{2}^{2}&\cdots&x_{r}^{2}\\
\vdots&\vdots&\vdots&\vdots\\
x_{1}^{r}&x_{2}^{r}&\cdots&x_{r}^{r}\\
\end{pmatrix}.
\end{equation}
It is easy to see that $\Delta_r$ is a highest weight vector for
both $gl(m)$ and $gl(n)$ and its weights are respectively
$\sum_{i=1}^r\epsilon_i$ and $\sum_{i=1}^r\tilde{\epsilon}_i$.
This weight corresponds to the Young diagram $$
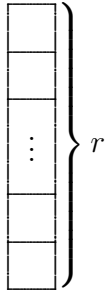
{\beginpicture
\setcoordinatesystem units <1.5pc,1.5pc> point at 0 0 \setplotarea
x from 0 to 1.5, y from 0 to 6 \plot 0 0 1 0 1 1 0 1 0 0 / \plot 1
1 1 2 0 2 0 1 / \plot 1 2 1 4 0 4 0 2 / \plot 1 4 1 5 0 5 0 4 /
\plot 1 5 1 6 0 6 0 5 / \put{$\vdots$} <0pt,2pt> at 0.5 3
\put{$\left.\hbox to 0pt{ \vrule width 0 pc height 4.4pc depth
4.4pc\hfil}
  \right\}r$} at 1.5 3
\endpicture}$$
That is, $\Delta_r$ is the highest weight vector for
$\Lambda^r(\C^m)\otimes\Lambda^r(\C^n)$ inside
$S(\C^m\otimes\C^n)$, the tensor product of the $r$-th fundamental
representations of $gl(m)$ and $gl(n)$.

Let $\lambda$ be a Young diagram as in \eqnref{standard} with
length not exceeding ${\rm min}(m,n)$. The set of highest weight
vectors in $S(C^m\otimes\C^n)$ form an abelian  semigroup, and the product
$\Delta_{\lambda'_1}\Delta_{\lambda'_2}\cdots\Delta_{\lambda'_t}$
is a highest weight vector for the irreducible representation in
$S(C^m\otimes\C^n)$ corresponding to the Young diagram $\lambda$.

On the other hand, the skew-symmetric algebra $\La (\C^m \otimes
\C^n)$ admits an induced $gl(m) \times gl (n)$ action. Following
Howe \cite{H2}, we have the multiplicity-free decomposition

\begin{equation*}
\La (\C^m\otimes\C^n) \cong \sum_{\la } V_m^{\lambda}\otimes
V_n^{\lambda '},
\end{equation*}
where the summation runs over Young diagrams $\lambda$ of length
not exceeding $m$ and of width not exceeding $n$.

Denote by $\eta_i^j, 1 \le i \le m, 1 \le j \le n$ the standard
basis for $\C^m\otimes\C^n$ in the consideration of skew-symmetric
algebra. The highest weight vector for the $gl(m) \times gl
(n)$-module $V_m^{\lambda}\otimes V_n^{\lambda '}$ inside
$\La(\C^m\otimes\C^n)$ is given by
\begin{eqnarray*}
 \eta_1^1 \eta_1^2 \cdots \eta_1^{\la_1}  \;
 \eta_2^1 \eta_2^2 \cdots \eta_2^{\la_2}  \ldots
 \eta_l^1 \eta_l^2 \cdots \eta_l^{\la_l}
\end{eqnarray*}
where $l $ is the length of $\la$.

Intimately related to the Howe duality is the Schur duality, which
we review below. Consider the standard representation of $GL(m)$
on $\C^m$. It induces an action on the $k$-th tensor power
$\otimes^k\C^m$. Now the symmetric group $S_k$ in $k$ letters acts
on $\otimes^k\C^m$ in a natural way.  These two actions commute
and we may thus decompose $\otimes^k\C^m$ into a direct sum of
irreducible $GL(m)\times S_k$-module.  Recalling that the
irreducible representations of symmetric group $S_k$ admit
parameterization by Young diagrams of weight $k$, Schur duality
states that
\begin{equation*}
\otimes^k\C^m\cong\sum_{\lambda } V_m^\lambda\otimes M_k^\lambda,
\end{equation*}
where the summation is over Young diagrams $\lambda$ of size $k$
and of length not exceeding $m$.  Here $M_k^\lambda$ is the
irreducible representation of $S_k$ corresponding to the Young
diagram $\lambda$.

Further well known examples of a multiplicity-free action of
$gl(m)$ that are of interest to us are as follows: consider the
action of $gl(m)$ on the symmetric square $S^2\C^m$ and
skew-symmetric square $\Lambda^2\C^m$.  We have an induced action
on their respective symmetric algebras $S(S^2\C^m)$ and
$S(\Lambda^2\C^m)$. Explicitly, the decomposition of these spaces
as $gl(m)$-modules is as follows (cf.~\cite{H2}, \cite{GW}):
\begin{align}
S(S^2\C^m)&\cong\sum_{l(\la)\le m}V^{2\la}_m,\label{symsym}\\
S(\Lambda^2\C^m)&\cong\sum_{l(\la)\le \frac{m}{2}}V^{(2\la)'}_m.\label{symskew}
\end{align}
Explicit formulas for the highest weight vectors in either cases
are well known (cf.~\cite{H2}) and are given in \remref{generic}.

One may also consider the decompositions of the skew-symmetric
algebra of $S^2\C^m$ and $\Lambda^2\C^m$.  In order to describe
the highest weights that appear in these decompositions we need a
few terminology.  The Young diagram associated to the partition
$\la=(k+1,1,,\ldots,1)$ of length $k\ge 1$ is called a {\it
$(k+1,k)$-hook}.  We will sometimes also call this $(k+1,k)$-hook
{\it a hook of shape $(k+1,k)$}. Assuming that $k>l$ we may form a
new Young diagram by ``nesting'' the $(l+1,l)$-hook inside the
$(k+1,k)$-hook.  The resulting partition of length $k$ is
$(k+1,l+2,2,\ldots,2,1,\dots,1)$, where $2$ appears $l-1$ times
and $1$ appears $k-l-1$ times.  Similarly a sequence of hooks of
shapes $(k_1+1,k_1),\ldots,(k_s+1,k_s)$ with $k_i>k_{i+1}$ for
$i=1,\ldots,s-1$ may be nested, and the resulting partition has
length $k_1$.  In consistency with the terminology used we call
the partition $(k,1,\ldots,1)$ of length $k+1$ {\it a hook of
shape $(k,k+1)$} or  a {\it $(k,k+1)$-hook}. Nesting of hooks of
shapes $(k_1,k_1+1),\ldots,(k_s,k_s+1)$ with $k_i>k_{i+1}$ for
$i=1,\ldots,s-1$ is done in an analogous fashion.

Now we can state the following multiplicity-free decompositions of
$gl(m)$-modules (cf.~\cite{H2}, \cite{GW}):
\begin{align}
\Lambda({S^2\C^m}) \cong \sum_{\lambda} V_m^{\lambda},\label{skewsym}\\
\Lambda({\Lambda^2\C^m}) \cong \sum_{\mu} V_m^{\mu},\label{skewskew}
\end{align}
where $\lambda$ (respectively $\mu$) is over all partitions with
$l(\la)\le m$ (respectively $l(\mu)\le m)$ such that $\la$
(respectively $\mu$) is obtained by nesting a sequence of
$(k+1,k)$-hooks (respectively of  $(k,k+1)$-hooks).

\section{Multiplicity-free Actions of the general linear Lie superalgebra}
\label{sect_superdual}

Let $\C^{m|n}$ denote the superspace of superdimension $m|n$.
Recall that this means that $\C^{m|n}$ is a $\Z_2$-graded space,
where the even subspace has dimension $m$ and the odd subspace has
dimension $n$.  The space of linear maps from $\C^{m|n}$ to itself
can be regarded as the space of $(m+n)\times(m+n)$ matrices with
an induced $\Z_2$-gradation, which gives it a natural structure as
a Lie superalgebra, denoted by $gl(m|n)$.  We have a triangular
decomposition $gl(m|n)=gl(m|n)_{-1}+ gl(m|n)_0+gl(m|n)_1$, where
$gl(m|n)_{\pm 1}$ denote the set of upper and lower triangular
matrices and $gl(m|n)_0$ denotes the set of diagonal matrices.
Given an $m+n$ tuple of complex numbers
$(a_1,\ldots,a_m;b_1,\ldots,b_n)$, we associate an irreducible
$gl(m|n)$-module $V_{m|n}$ of highest weight
$(a_1,\ldots,a_m;b_1,\ldots,b_n)$ (with respect to the standard
Borel subalgebra $gl(m|n)_0+gl(m|n)_{1}$). It is well known
(cf.~e.g.~\cite{K}) that the module $V_{m|n}$ is
finite-dimensional if and only if $(a_1,\ldots,a_m;
b_1,\ldots,b_n)$ satisfies the conditions
$a_i-a_{i+1},b_j-b_{j+1}\in\Z_+$, for all $i=1,\ldots,m-1$ and
$j=1,\ldots,n-1$.

Let $\C^{p|q}$ and $\C^{m|n}$ denote complex superspaces of
superdimensions $p|q$ and $m|n$, respectively.  We will now
describe a duality between the Lie superalgebras $gl(p|q)$ and
$gl(m|n)$. Our starting point is Schur duality for Lie
superalgebra $gl(m|n)$.

Schur duality for the Lie superalgebra $gl(m|n)$ was studied in
\cite{Se}.  Below we will recall the main result for the
convenience of the reader.  Let $\C^{m|n}$ denote the standard
$gl(m|n)$-module.  We may, as in the classical case, consider the
$k$-th tensor power $\otimes^k\C^{m|n}$ which admits a natural
action of the Lie superalgebra $gl(m|n)$. On the other hand the
symmetric group $S_k$ acts naturally on $\otimes^k\C^{m|n}$ by
permutations with appropriate signs (corresponding to the
permutations of odd elements in $\C^{m|n}$). It is easy to check
that the actions of $gl(m|n)$ and $S_k$ commute with each other,
cf. Sergeev \cite{Se} (also see Berele-Regev \cite{BR} for a more
detailed study).

\begin{thm} [Sergeev] \label{Sergeev}
As a $gl(m|n)\times S_k$-module we have
\begin{equation*}
\otimes^k \C^{m|n} \cong \sum_{\lambda} V^{\lambda}_{m|n} \otimes
M_k^{\lambda},
\end{equation*}
where $\lambda$ is summed over Young diagrams of size $k$ such
that $\lambda_{m+1}\le n$, $M_k^{\lambda}$ is the irreducible
$S_k$-module parameterized by $\lambda$, and $V^{\lambda}_{m|n}$
denotes the irreducible $gl(m|n)$-module with highest weight
$(\lambda_1, \lambda_2, \ldots, \lambda_m; \lef\lambda'_1-m \ri,
\ldots, \lef \lambda'_n-m \ri)$, where we denote $\lef l\ri =l$
for $l\in\Z_+$ and $\lef l \ri=0$, otherwise.
\end{thm}

The symmetric algebra $S(\C^{p|q}\otimes\C^{m|n})$ is by
definition equal to the tensor product of the symmetric algebra of
the even part of $\C^{p|q}\otimes\C^{m|n}$ and the skew-symmetric
algebra of the odd part of $\C^{p|q}\otimes\C^{m|n}$. It admits a
natural gradation
\begin{eqnarray*}
S( \C^{p|q} \otimes\C^{m|n}) = \sum_{k \geq 0} S^k
(\C^{p|q}\otimes\C^{m|n})
\end{eqnarray*}
by letting the degree of the basis elements of
$\C^{p|q}\otimes\C^{m|n}$ be $1$. The natural actions of $gl(p|q)$
on $\C^{p|q}$ and $gl(m|n)$ on $\C^{m|n}$ induce commuting actions
on the $k$-th symmetric algebra $S^k(\C^{p|q} \otimes \C^{m|n})$.
Indeed $gl(p|q)$ and $gl(m|n)$ are mutual centralizers in
$gl(\C^{p|q} \otimes \C^{m|n})$. We obtain the following theorem
by an analogous argument as in \cite{H2}.

\begin{thm}\label{supergl-duality}
The symmetric algebra $S(\C^{p|q}\otimes\C^{m|n})$ is
multiplicity-free as a module over $gl(p|q)\times gl(m|n)$. More
explicitly, we have the following decomposition
\begin{equation*}
S(\C^{p|q}\otimes\C^{m|n}) \cong \sum_{\lambda}
V^{\lambda}_{p|q}\otimes V^{\lambda}_{m|n},
\end{equation*}
where the sum is over Young diagrams $\lambda$ satisfying
$\lambda_{p+1}\le q$ and $\lambda_{m+1}\le n$.  Here the highest
weight of the module $V^{\lambda}_{p|q}$ (respectively
$V^{\lambda}_{m|n}$) is given by $(\lambda_1, \ldots, \lambda_p;
\lef \lambda_1'-p \ri, \ldots, \lef \lambda'_q-p \ri)$ (resp.
$(\lambda_1, \ldots, \lambda_m; \lef \lambda_1'-m \ri, \ldots,
\lef \lambda'_n-m \ri)$).
\end{thm}

\begin{proof}
By the definition of the $k$-th supersymmetric algebra we have
\begin{equation*}
S^k (\C^{m|n}\otimes\C^{p|q}) \cong ((\otimes^k\C^{m|n}) \otimes
(\otimes^k\C^{p|q}))^{\Delta_k},
\end{equation*}
where $\Delta_k$ is the diagonal subgroup of $S_k\times S_k$.  By
\thmref{Sergeev} we have therefore

\begin{align*}
S(\C^{m|n}\otimes\C^{p|q})&\cong
\sum_{k=0}^{\infty}((\sum_{|\lambda|=k}V_{m|n}^{\lambda}\otimes
M_k^{\lambda})\otimes(\sum_{|\mu|=k} V_{p|q}^{\mu}\otimes
M_k^{\mu}))^{\Delta_k}\\ &\cong \sum_{k=0}^{\infty}
\sum_{|\lambda|=|\mu|=k}(V_{m|n}^{\lambda}\otimes
V_{p|q}^\mu)\otimes (M_k^{\lambda}\otimes M_k^{\mu})^{\Delta_k}\\
&\cong\sum_{k=0}^\infty\sum_{|\lambda|=k}
(V_{m|n}^{\lambda}\otimes V_{p|q}^\lambda) \\ &\cong\sum_{\la}
(V_{m|n}^{\lambda}\otimes V_{p|q}^\lambda),
\end{align*}
where $\lambda$ in the previous line is summed over all Young
diagrams satisfying the conditions $\lambda_{m +1} \leq n$ and
$\la_{p+1}\le q$. The second to last equality follows from the
well-known fact that $M_k^{\lambda}$ is a self-contragredient
module.
\end{proof}

\begin{rem}
\begin{enumerate}
\item This theorem (except the explicit formula for
the highest weights) was first obtained in Brini et al \cite{BPT}
in a combinatorial approach. It is also obtained independently
recently by Sergeev.
\item When $n =q =0$, we recover the $(gl(p),gl(m))$-duality in the
symmetric algebra case. When $q =m = 0$ we recover the
$(gl(p),gl(n))$-duality in the skew-symmetric algebra case.
\item One can easily show that the $(gl(m|n), gl (k) )$-duality
implies the $(gl(m|n), S_k )$ Schur duality (\thmref{Sergeev}),
using an argument of Howe (cf. 2.4, \cite{H2}).
\end{enumerate}
\end{rem}

The next corollary is immediate from \thmref{supergl-duality}.

\begin{cor}
The image of the action of the universal enveloping algebras of
$gl(p|q)$ and $gl(m|n)$ on $S^k(\C^{p|q}\otimes\C^{m|n})$ are
double commutants.
\end{cor}

\begin{thm}\label{skewsuperduality}
The skew-symmetric algebra $\Lambda(\C^{p|q}\otimes\C^{m|n})$ is
multiplicity-free as a module over $gl(p|q)\times gl(m|n)$. More
explicitly, we have the following decomposition
\begin{equation*}
\Lambda(\C^{p|q}\otimes\C^{m|n}) \cong \sum_{\lambda}
V^{\lambda}_{p|q}\otimes V^{\lambda'}_{m|n},
\end{equation*}
where the sum is over Young diagrams $\lambda$ satisfying
$\lambda_{p+1}\le q$ and $\lambda'_{m+1}\le n$.  Here the highest
weight of the module $V^{\lambda}_{p|q}$ (respectively
$V^{\lambda'}_{m|n}$) is given by $(\lambda_1, \ldots, \lambda_p;
\lef \lambda_1'-p \ri, \ldots, \lef \lambda'_q-p \ri)$ (resp.
$(\lambda'_1, \ldots, \lambda'_m; \lef \lambda_1-m \ri, \ldots,
\lef \lambda_n-m \ri)$).
\end{thm}

\begin{proof}
By the definition of the $k$-th skew-symmetric algebra we have
\begin{equation*}
\Lambda^k (\C^{p|q}\otimes\C^{m|n}) \cong ((\otimes^k\C^{p|q}) \otimes
(\otimes^k\C^{m|n}))^{\Delta_k\sim},
\end{equation*}
where $\Delta_k$ is  the diagonal subgroup of $S_k\times S_k$ and
$(\otimes^k\C^{m|n}))^{\Delta_k\sim}$ is the subspace of
$(\otimes^k\C^{m|n}))$ that transforms according to the sign
character of $\Delta_k$. By \thmref{Sergeev} we have therefore

\begin{align*}
\Lambda(\C^{p|q}\otimes\C^{m|n})&\cong \sum_{k=0}^{\infty} \left(
(\sum_{|\lambda|=k}V_{p|q}^{\lambda}\otimes
M_k^{\lambda})\otimes(\sum_{|\mu|=k} V_{m|n}^{\mu}\otimes
M_k^{\mu}) \right)^{\Delta_k\sim}\\ &\cong \sum_{k=0}^{\infty}
\sum_{|\lambda|=|\mu|=k}(V_{p|q}^{\lambda}\otimes
V_{m|n}^\mu)\otimes (M_k^{\lambda}\otimes
M_k^{\mu})^{\Delta_k\sim}\\
&\cong\sum_{k=0}^\infty\sum_{|\lambda|=k} V_{p|q}^{\lambda}\otimes
V_{m|n}^{\lambda'} \\ &\cong\sum_{\la} V_{p|q}^{\lambda}\otimes
V_{m|n}^{\lambda'},
\end{align*}
where $\lambda$ in the previous line is summed over all Young
diagrams satisfying the conditions $\lambda_{p+1} \leq q$ and
$\la'_{m+1}\le n$. The second to last equality follows from the
following well known facts: $M_k^{\lambda}$ is a self-contragredient
module and tensoring the module $M_k^{\lambda}$ with the sign
character yields the module $M_k^{\lambda'}$.
\end{proof}

\begin{rem}
Of course it follows from \thmref{skewsuperduality} that the image
of the action of the universal enveloping algebras of $gl(p|q)$
and $gl(m|n)$ on $\Lambda^k(\C^{p|q}\otimes\C^{m|n})$ are also
double commutants.
\end{rem}

The following corollary turns out to be very useful later on in
order to check that a given vector is indeed a highest weight
vector inside $S(\C^{p|q}\otimes\C^{m|n})$.

\begin{cor}  \label{cor_hwt}
Assume a vector $v \in S(\C^{p|q}\otimes\C^{m|n})$ has the weight
$\la$ with respect to $gl(p|q)\times gl(m|n)$ associated to a
Young diagram $\la$ satisfying $\lambda_{p+1}\le q$ and
$\lambda_{m+1}\le n$. If $v$ is a highest weight vector for
$gl(p|q)$, then it is for $gl(m|n)$ as well.
\end{cor}

\begin{proof}
Since $v$ is a highest weight vector for $gl(p|q)$ with weight
$\la$, it belongs to the subspace $W$ of
$S(\C^{p|q}\otimes\C^{m|n})$ which consists of vectors with weight
$\la$ annihilated by the standard Borel in $gl(p|q)$. By
Theorem~\ref{supergl-duality}, $W$ is isomorphic to
$V_{m|n}^{\la}$ as a $gl(m|n)$-module. There exists a unique
vector (up to scalar multiple) in $V_{m|n}^{\la}$ which has weight
$\la$, which is the highest weight vector. By assumption $v$ has
weight $\la$ as a $gl(m|n)$-module, so it is a highest weight
vector for $gl(m|n)$.
\end{proof}

The description of highest weight vectors of the irreducible
$gl(p|q)\times gl(m|n)$-modules in the symmetric algebra
turns out to be much more subtle than in the classical Howe
duality case and we will deal with this question in \secref{jointhwv}.

Next consider the symmetric square $S^2\C^{m|n}$ of the natural
representation of $gl(m|n)$.  The following theorem can be proved
by an analogous argument as in \cite{H2}. This result was also
obtained independently recently by Sergeev. We omit the proof
since it is in any case parallel to the proof of
\thmref{skew-symmetric} below.

\begin{thm}\label{symmetric}
The symmetric algebra of the symmetric square of the
natural representation $\C^{m|n}$ of the Lie superalgebra
$gl(m|n)$ is a completely reducible multiplicity-free
$gl(m|n)$-module.  More precisely we have the following
decomposition
\begin{equation*}
S^k(S^2\C^{m|n})=\sum_{\la}V_{m|n}^\la,
\end{equation*}
where the summation is over all partitions $\la$ into even parts
of size $2k$ and $\la_{m+1}\le n$.
\end{thm}

Now $S^2\C^{m|n}$ reduces to $S^2\C^m$ in the case when $n=0$, and
to $\Lambda^2\C^n$ in the case when $m=0$, the symmetric and
skew-symmetric square of the natural representation of $gl(m)$ and
$gl(n)$, respectively. Thus one obtains as a corollary the
classical multiplicity-free decompositions of their respective
symmetric algebras, namely \eqnref{symsym} and \eqnref{symskew}.
Again the question of obtaining explicit formulas for the highest
weight vectors inside $S(S^2\C^{m|n})$ is substantially more
subtle than in the non-super case.  We will give these in
\secref{hwvector}.

\begin{thm}\label{skew-symmetric}
The skew-symmetric algebra of the symmetric square of the natural
representation $\C^{m|n}$ of the Lie superalgebra $gl(m|n)$ is a
completely reducible multiplicity-free $gl(m|n)$-module.  More
precisely we have the following decomposition
\begin{equation*}
\Lambda^k({S^2\C^{m|n}})=\sum_{\la}V_{m|n}^\la,
\end{equation*}
where the summation is over all partitions $\la$ of size $2k$,
which are obtained by nesting $(l+1,l)$-hooks with $\la_{m+1}\le
n$.
\end{thm}

\begin{proof}
Our argument follows closely the one given in the proof of Theorem
4.4.2 in \cite{H2} with \thmref{Sergeev} replacing the classical
Schur duality.  Let $D_k$ denote the the subgroup of $S_{2k}$,
which preserves the partition
$\{\{1,2\},\{3,4\},\ldots,\{2k-1,2k\}\}$ of $2k$. Note that $D_k$
is isomorphic to a semidirect product of $S_k$ and $(\Z_2)^k$,
where $\Z_2$ acts by interchanging $2j-1$ with $2j$ and $S_k$ acts
by permuting the pairs. Let ${\rm sign\sim}$ denote the character
on $D_k$ which is trivial on $(\Z_2)^k$, but transforms by the
sign character on $S_k$.  We observe that
\begin{equation*}
\Lambda^{2k}(S^2\C^{m|n})\cong(\bigotimes^{2k}\C^{m|n})^{D_k,\rm
sign\sim}.
\end{equation*}
Thus using \thmref{Sergeev} we obtain
\begin{equation*}
\Lambda^{2k}(S^2\C^{m|n}) \cong\sum_{|\la|=2k}(V^\la_{m|n}\otimes
M^\la_{2k})^{D_k,\rm sign\sim}\cong
\sum_{|\la|=2k}V^\la_{m|n}\otimes (M^\la_{2k})^{D_k,\rm sign\sim}.
\end{equation*}
Now by Theorem A1.4 of \cite{H2} the space $(M^\la_{2k})^{D_k,\rm
sign\sim}$ is non-zero if and only if $\la$ is constructed from
nesting hooks of types $(l+1,l)$, in which case it is
one-dimensional.
\end{proof}

Similarly we obtain as a corollary the classical multiplicity-free
decompositions \eqnref{skewsym} and \eqnref{skewskew}.

\begin{rem}
The character of $V^\la_{m|n}$ is defined as the trace of the
action of the diagonal matrix ${\rm
diag}(x_1,\ldots,x_m;y_1,\ldots,y_n)$ in $gl(m|n)$ on
${V^\la_{m|n}}$ and according to \cite{BR} is given by so-called
hook Schur functions ${\rm HS}_{\la}(x,y)$ (see \cite{BR} for
definition). Thus, comparing the characters of both sides of
\thmref{supergl-duality} and \thmref{skewsuperduality},
respectively, with $x=(x_1,\ldots,x_p)$, $y=(y_1,\ldots,y_q)$,
$u=(u_1,\ldots,u_m)$ and $v=(v_1,\ldots,v_n)$ we obtain the
following combinatorial identities:
\begin{align*}
\sum_{\la}{\rm HS}_{\la}(x,y){\rm HS}_{\la}(u,v)&=
\prod_{i,j,k,l}(1-x_iu_{k})^{-1}(1-y_jv_{l})^{-1}(1+x_iv_l)
(1+y_ju_k),\\
\sum_{\la}{\rm HS}_{\la}(x,y){\rm HS}_{\la'}(u,v)&=
\prod_{i,j,k,l}(1+x_iu_{k})(1+y_jv_{l})(1-x_iv_l)^{-1}
(1-y_ju_k)^{-1},
\end{align*}
where $1\le i\le p$, $1\le j\le q$, $1\le k\le m$ and $1\le l\le
n$ with summation in the first identity over $\la$ such that
$\la_{p+1}\le q$ and $\la_{m+1}\le n$ and in the second one over
$\la$ such that $\la_{p+1}\le q$ and $\la'_{m+1}\le n$.
Now putting $y=v=0$ in the first identity we obtain the
classical Cauchy identity, while putting respectively
$y=u=0$ the dual Cauchy identity (see e.g.~\cite{M}):
\begin{align*}
&\sum_{\la}s_\la(x)s_\la(y)=\prod_{i,j}(1-x_iy_j)^{-1},\\
&\sum_{\mu}s_{\mu}(x)s_{\la'}(y)=\prod_{i,j}(1+x_iy_j).
\end{align*}
\end{rem}

\begin{rem}
Similarly \thmref{symmetric} and \thmref{skew-symmetric}
give rise to the following combinatorial identities
($x=(x_1,\ldots,x_m)$ and $y=(y_1,\ldots,y_n)$):
\begin{align*}
\sum_{\la}{\rm HS}_{\la}(x,y)&=\prod_{i\le
i',j<j'}(1-x_ix_{i'})^{-1}
(1-y_jy_{j'})^{-1}\prod_{i,j}(1+x_iy_j),\\
\sum_{\mu}{\rm HS}_{\mu}(x,y)&=\prod_{i\le i',j<j'}(1+x_ix_{i'})
(1+y_jy_{j'})\prod_{i,j}(1-x_iy_j)^{-1},
\end{align*}
where in the first identity the sum is over all partitions $\la$
with even rows such that $\la_{m+1}\le n$ and in the second over
all partitions $\mu$ that can be obtained by nesting
$(k+1,k)$-hooks such that $\mu_{m+1}\le n$ and $1\le i,i'\le m$,
$1\le j,j'\le n$. Putting either $x=0$ or $y=0$ in these two
identities we obtain the following classical Schur function
identities (see e.g.~\cite{M}), which correspond to the
decompositions in \eqnref{symsym}, \eqnref{symskew},
\eqnref{skewsym} and \eqnref{skewskew}, respectively:
\begin{align*}
&\sum_{l(\la)\le m}s_{2\la}=\prod_{1\le i\le m
}(1-x_i^2)^{-1}\prod_{1\le i<j\le m }(1-x_ix_j)^{-1},\\
&\sum_{l(\mu')\le\frac{m}{2}}s_{(2\mu)'}=\prod_{1\le i<j\le
m}(1-x_ix_j)^{-1},\\ &\sum_{\rho}s_{\rho}=\prod_{1\le i\le m
}(1+x_i^2)\prod_{1\le i<j\le m }(1+x_ix_j),\\
&\sum_{\pi}s_{\pi}=\prod_{1\le i<j\le m }(1+x_ix_j),\\
\end{align*}
where $\rho$ (respectively $\pi$) above is summed over all nested
sequences of hooks of shape $(k+1,k)$ with $k\le m$ (respectively
of hooks of shape $(k,k+1)$ with $k\le m-1$).
\end{rem}

\section{Construction of highest weight vectors in
$S(\C^{p|q}\otimes\C^{m|n})$}\label{jointhwv}

This section is devoted to the construction of the
highest weight vectors of $gl(p|q)\times gl(m|n)$ inside the
symmetric algebra of $\C^{p|q}\otimes\C^{m|n}$.  We will
divide this section into several cases.  Before we embark on this
task we will set the notation to be used throughout this section.

We let $e^1,\ldots,e^p;f^1,\ldots,f^q$ denote the standard
homogeneous basis for the standard $gl(p|q)$-module.  Here $e^i$
are even, while $f^j$ are odd basis elements.  Similarly we let
$e_1,\ldots,e_m;f_1,\ldots,f_n$ denote the standard homogeneous
basis for the standard $gl(m|n)$-module.  The weights of $e^i,
f^j,e_l$ and $f_k$ are denoted by $\tilde{\epsilon}_i$,
$\tilde{\delta}_j$, $\epsilon_l$ and $\delta_k$, for $1\le i\le
p$, $1\le j\le q$, $1\le l\le m$ and $1\le k\le n$, respectively.
We set
\begin{equation}\label{generators}
x_l^i:=e_l\otimes e^i;\ \xi_l^j:=e_l\otimes f^j;\
\eta_k^i:=f_k\otimes e^i;\ y_k^j:=f_k\otimes f^j.
\end{equation}
We will denote by $\C[\x,\xibf,\etabf,\y]$ the polynomial
superalgebra generated by \eqnref{generators}. The commuting pair
of $gl(p|q)$ and $gl(m|n)$ may be realized as first order
differential operators as follows ($1\le i,i'\le p;1\le l,l'\le q$
and $1\le s,s'\le m;1\le k,k'\le n$):
\begin{align}
&\sum_{j=1}^m x_{j}^{i}\frac{\partial}{\partial
x_j^{i'}}+\sum_{j=1}^n\eta_j^{i}\frac{\partial}{\partial\eta_j^{i'}},
\quad \sum_{j=1}^m \xi_{j}^{l}\frac{\partial}{\partial
\xi_j^{l'}}+\sum_{j=1}^n y_j^{l}\frac{\partial}{\partial
y_j^{l'}},\label{glpq}\\
 &\sum_{j=1}^m
x_{j}^{i}\frac{\partial}{\partial
\xi_j^{l'}}+\sum_{j=1}^n\eta_j^{i}\frac{\partial}{\partial
y_j^{l'}}, \quad \sum_{j=1}^m \xi_{j}^{l}\frac{\partial}{\partial
x_j^{i'}}+\sum_{j=1}^n
y_j^{l}\frac{\partial}{\partial\eta_j^{i'}},\nonumber\\
&\sum_{j=1}^p x_{s}^{j}\frac{\partial}{\partial
x_{s'}^j}+\sum_{j=1}^q \xi_{s}^{j}\frac{\partial}{\partial
\xi_{s'}^j},\quad
\sum_{j=1}^p\eta_{k'}^{j}\frac{\partial}{\partial\eta_k^j}+
\sum_{j=1}^q y_{k'}^{j}\frac{\partial}{\partial
y_k^j},\label{glmn}\\
 &\sum_{j=1}^p
x_s^j\frac{\partial}{\partial\eta_k^j}- \sum_{j=1}^q
\xi_s^j\frac{\partial}{\partial y_k^j},\quad
\sum_{j=1}^p\eta_k^j\frac{\partial}{\partial x_s^j}- \sum_{j=1}^q
y_k^j\frac{\partial}{\partial \xi_s^j}.\nonumber
\end{align}
\eqnref{glpq} spans a copy of $gl(p|q)$, while \eqnref{glmn} spans a copy of
$gl(m|n)$.

Our Cartan subalgebras of $gl(p|q)$ and $gl(m|n)$ are spanned, respectively, by
\begin{equation*}
\sum_{j=1}^m x_{j}^{i}\frac{\partial}{\partial
x_j^{i}}+\sum_{j=1}^n\eta_j^{i}\frac{\partial}{\partial\eta_j^{i}},
\sum_{j=1}^m \xi_{j}^{l}\frac{\partial}{\partial \xi_j^{l}}+\sum_{j=1}^n
y_j^{l}\frac{\partial}{\partial y_j^{l}},
\end{equation*}
and
\begin{equation*}
\sum_{j=1}^p x_{s}^{j}\frac{\partial}{\partial x_{s}^j}+\sum_{j=1}^q
\xi_{s}^{j}\frac{\partial}{\partial \xi_{s}^j},\
\sum_{j=1}^p\eta_{k}^{j}\frac{\partial}{\partial\eta_k^j}+
\sum_{j=1}^q y_{k}^{j}\frac{\partial}{\partial y_k^j},
\end{equation*}
while the nilpotent radicals are respectively generated by the
simple root vectors

\begin{eqnarray}\label{radicalpq}
\sum_{j=1}^m x_{j}^{i-1}\frac{\partial}{\partial
x_j^{i}}+\sum_{j=1}^n\eta_j^{i-1}\frac{\partial}{\partial\eta_j^{i}},
\sum_{j=1}^m \xi_{j}^{l-1}\frac{\partial}{\partial
\xi_j^{l}}+\sum_{j=1}^n y_j^{l-1}\frac{\partial}{\partial
y_j^{l}}, \nonumber   \\
 \sum_{j=1}^m x_{j}^{p}\frac{\partial}{\partial
\xi_j^{1}}+\sum_{j=1}^n\eta_j^{p}\frac{\partial}{\partial
y_j^{1}}, \quad 1< i \le p, 1< l \le q.
\end{eqnarray}
and
\begin{align}\label{radicalmn}
\sum_{j=1}^p x_{s-1}^{j}\frac{\partial}{\partial
x_{s}^j}+\sum_{j=1}^q \xi_{s-1}^{j}\frac{\partial}{\partial
\xi_{s}^j},\
\sum_{j=1}^p\eta_{k-1}^{j}\frac{\partial}{\partial\eta_k^j}+
\sum_{j=1}^q y_{k-1}^{j}\frac{\partial}{\partial y_k^j},
 \nonumber  \\
 \sum_{j=1}^p x_m^j\frac{\partial}{\partial\eta_1^j}- \sum_{j=1}^q
\xi_m^j\frac{\partial}{\partial y_1^j},\quad 1< s \le m,1< k \le
n.
\end{align}
With these conventions, we may thus identify
$S(\C^{p|q}\otimes\C^{m|n})$ with the polynomial superalgebra
$\C[\x,\xibf,\etabf,\y]$ (as $gl(p|q)\times gl(m|n)$-modules).

\subsection{Highest Weight Vectors: the Case $q=0$} \label{sect_hwto}

In this section we will describe the highest weight vectors for
$gl(p) \times gl(m|n)$ in the symmetric algebra
$S(\C^p\otimes\C^{m|n})$, i.e. $q=0$ case. The space
$S(\C^p\otimes\C^{m|n})$ is identified with $\C[\x,\etabf]$, and
\eqnref{radicalpq} and \eqnref{radicalmn} reduce to
\begin{eqnarray}
 &\sum_{j=1}^m x_{j}^{i-1} \frac{\partial}{\partial x_j^{i}}
+\sum_{j=1}^n\eta_j^{i-1}\frac{\partial}{\partial\eta_j^{i}},
\label{radicalp}\\
 &\sum_{j=1}^p
x_{s-1}^{j}\frac{\partial}{\partial x_{s}^j},\
\sum_{j=1}^p\eta_{k-1}^{j}\frac{\partial}{\partial\eta_k^j},\
\sum_{j=1}^p
x_m^j\frac{\partial}{\partial\eta_1^j}\label{radicalpmn},
\end{eqnarray}
respectively.  Now by \thmref{supergl-duality} a highest weight
representation $V_p^{\la}\otimes V_{m|n}^{\la}$ of $gl(p)\times
gl(m|n)$ appears in the decomposition of
$S^{k}(\C^p\otimes\C^{m|n})$ if and only if $\la$ is of size $k$
and of length at most $p$ such that $\la_{m+1}\le n$.

We will consider two cases separately, namely $m\ge p$ and $m<p$.

We begin with the case of $m\ge p$.  Here the condition
$\la_{m+1}\le n$ is an empty condition.  So we are looking for
homogeneous polynomials of degree $k$ in $\C[\x,\etabf]$,
annihilated by all vectors of \eqnref{radicalp} and
\eqnref{radicalpmn}, and having $gl(p)$- and $gl(m|n)$-weight
$\la$ of length not exceeding $p$. If $\la$ is such a weight, then
$\la_i'\le p$, where we recall that $\la ' = (\la_1' , \ldots,
\la_t')$ denotes the transpose of $\la$. It is easy to see that
the product $\Delta_{\la'_1}\cdots\Delta_{\la'_t}$ is annihilated
by all vectors of \eqnref{radicalp} and \eqnref{radicalpmn}, where
we recall that $\Delta_r$ is defined in \eqnref{deltar}. It is
straightforward to check that its weight is exactly $\la$.

\begin{thm}
In the case when $m\ge p$, all $gl(p) \times gl(m|n)$ highest
weight vectors in $ \C[\x, \etabf]$ form an abelian semigroup generated by
$\Delta_r$, for $r=1,\ldots,p$. The highest weight vector
associated to the weight $\la$ is given by the product
$\Delta_{\la'_1}\cdots\Delta_{\la'_t}$.
\end{thm}

We now consider the case $p>m$.  In this case the condition
$\la_{m+1}\le n$ is no longer an empty condition.  Obviously the
highest weight vectors associated to Young diagrams $\la$ with
$\la_{m+1}=0$ can be obtained just as in the previous case.

Now suppose $\la$ is a diagram of length exceeding $m$.  Let
$\lambda_1',\lambda_2',\ldots,\la'_{t}$ denote its column lengths
as usual.  We have $p \ge \la'_1 \ge \la'_2 \ldots \ge \la'_t$ and
$m\ge\la'_{n+1}$.  For $m\le r\le p$, the following determinant of
an $r\times r$ matrix plays a fundamental role in this paper:
\begin{equation} \label{eq_det}
\Delta_{k,r}:={\rm det}
\begin{pmatrix}
x_1^1&x_1^2&\cdots &x_1^r\\
x_2^1&x_2^2&\cdots &x_2^r\\
\vdots&\vdots&\cdots &\vdots\\
x_m^1&x_m^2&\cdots &x_m^r\\
\eta_k^1&\eta_k^2&\cdots &\eta_k^r\\
\eta_k^1&\eta_k^2&\cdots &\eta_k^r\\
\vdots&\vdots&\cdots &\vdots\\
\eta_k^1&\eta_k^2&\cdots &\eta_k^r\\
\end{pmatrix},\quad k=1,\ldots,n.
\end{equation}
That is, the first $m$ rows are filled by the vectors $(x_j^1,
\ldots, x_j^r)$, for $j=1, \ldots, m$, in increasing order and the
last $r-m$ rows are filled with the same vector
$(\eta_k^1,\ldots,\eta_k^r)$. Since the matrix entries involve
Grassmann variables $\eta_k^i$, we must specify what we mean
by the determinant. By the determinant of a matrix

$$
 A :=\begin{pmatrix} a_1^1 & a_1^2 &\cdots&a_1^r\\
a_2^1&a_2^2&\cdots&a_2^r\\
\vdots&\vdots&\cdots &\vdots\\
a_r^1&a_p^2&\cdots&a_r^r\\
\end{pmatrix},
$$
whose matrix entries involve Grassmann variables $\eta_k^i$, we
will always mean the expression $\sum_{\sigma\in
S_r}(-1)^{p(\sigma)}a_1^{\sigma(1) }a_2^{\sigma (2) }\cdots
a_r^{\sigma (r) }$, where $p(\sigma)$ is the length of $\sigma$ in
the symmetric group $S_r$. In general it is not true that ${\rm
det}A={\rm det}A^t$.

\begin{rem}
The determinant (\ref{eq_det}) is always nonzero. It reduces to
(\ref{deltar}) when $m =r$, and reduces to (up to a scalar
multiple) $\eta_k^1 \cdots \eta_k^r$ when $m=0$.
\end{rem}

Now let $\la$ be a diagram of length at most $p$ such that
$\la_{m+1}\le n$.  It is thus of the following shape:
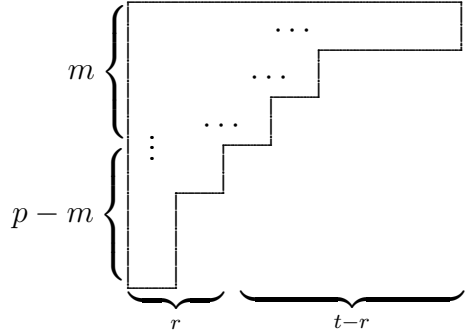
\begin{equation*}
{\beginpicture
\setcoordinatesystem units <1.5pc,1.5pc>
point at 0 3
\setplotarea x from 0 to 1.5, y from 0 to 6
\plot 0 0 1 0 1 2 2 2 2 3 3 3 3 4 4 4 4 5 7 5 7 6 0 6 0 0 /
\put{$\vdots$} <0pt,2pt> at 0.5 3
\put{$\cdots$} <0pt,2pt> at 3.5 5.3
\put{$\cdots$} <0pt,2pt> at 3 4.3
\put{$\cdots$} <0pt,2pt> at 2 3.3
\put{$\vdots$} <0pt,2pt> at 0.5 3
\put{$m\left\{\hbox to 0pt{
\vrule width 0 pc height 2.0pc depth 2.0pc\hfil}
  \right.$} at -.6 4.5
\put{$p-m\left\{\hbox to 0pt{
\vrule width 0 pc height 2pc depth 2pc\hfil}
  \right.$} at -1.2 1.5
\put{$\underbrace{\hbox to 3pc{\ \hfill}}_{r}$} at 1 -0.5
\put{$\underbrace{\hbox to 7pc{\ \hfill}}_{t-r}$} at 4.7 -0.5
\endpicture}
\end{equation*}
where $r$ is defined by $\la'_r>m$ and $\la'_{r+1}\le m$.
We can divide such a diagram into two diagrams, namely
\begin{equation}\label{splitq=0}
{\beginpicture
\setcoordinatesystem units <1.5pc,1.5pc> point at 0 3
\setplotarea x from 0 to 1.5, y from 0 to 6
\plot 0 0 1 0 1 2 2 2 2 6 0 6 0 0 /
\put{$m\left\{\hbox to 0pt{
\vrule width 0 pc height 2.0pc depth 2.0pc\hfil}
  \right.$} at -.6 4.5
\put{$p-m\left\{\hbox to 0pt{
\vrule width 0 pc height 2pc depth 2pc\hfil}
  \right.$} at -1.2 1.5
\put{$\underbrace{\hbox to 3pc{\ \hfill}}_{r}$} at 1 -0.5
\plot 7 3 8 3 8 4 9 4 9 5 12 5 12 6 7 6 7 3 /
\put{$m\left\{\hbox to 0pt{
\vrule width 0 pc height 2.0pc depth 2.0pc\hfil}
  \right.$} at 6.4 4.5
\put{$\underbrace{\hbox to 7pc{\ \hfill}}_{t-r}$} at 9.7 -0.5
\endpicture}
\end{equation}
Now the second diagram in \eqnref{splitq=0} has length not
exceeding $m$, so its associated highest weight vector is given by
the product $\Delta_{\la'_{r+1}}\cdots\Delta_{\la'_{t}}$. A
formula for the highest weight vector associated to the first
diagram in \eqnref{splitq=0} is given by the following
proposition. We will denote by $\prod_{k=1}^r\Delta_{k,\la'_k}$
the (ordered) product $\Delta_{1, \la'_1} \Delta_{2,\la'_2} \ldots
\Delta_{r,\la'_r}$.

\begin{prop}
Let $p\ge \la'_1\ge \la'_2\ge \ldots\ge \la'_r> m$.  Then
$\prod_{k=1}^r\Delta_{k,\la'_k}$ is a highest weight vector
associated to the first Young diagram in \eqnref{splitq=0}.
\end{prop}
\begin{proof}
Observe that $\prod_{k=1}^r\Delta_{k,\la'_k}$ has the same
$gl(p)\times gl(m|n)$-weight as the first Young diagram of
\eqnref{splitq=0}. Clearly $\prod_{k=1}^r\Delta_{k,\la'_k}$ is
non-zero. It is straightforward to verify that
$\prod_{k=1}^r\Delta_{k,\la'_k}$ is annihilated by the operators
in \eqnref{radicalp}. It follows from Corollary~\ref{cor_hwt} that
$\prod_{k=1}^r\Delta_{k,\la'_k}$ is also a highest weight vector
for $gl(m|n)$.
\end{proof}

Our next theorem follows by observing that the product of the
highest weight vectors corresponding to the two Young diagrams in
\eqnref{splitq=0} is non-zero and is a highest weight vector
associated to the Young diagram $\la$.

\begin{thm}\label{glpmn-duality}
Suppose that $m<p$.  An irreducible highest weight module
$V_{p}^\la\otimes V_{m|n}^\la$ appearing in $ \C[ \x,\eta]$ if and
only if $\la$ corresponds to a Young diagram $\la$ of length not
exceeding p and $\la_{m+1}\le n$.  Furthermore a highest weight
vector associated to such a $\la$ is given by
$$\prod_{k=1}^r\Delta_{k,\la'_k}\prod_{j=r+1}^t\Delta_{\la'_j},$$
where $r$ is defined by $\la'_r>m$ and $\la'_{r+1}\le m$.
\end{thm}

As a corollary we obtain the following useful combinatorial
identity, which will play an important role later on.

\begin{cor}\label{identity}
Let $x^i_l$ be even variables for $i=1,\ldots,p$ and
$l=1,\ldots,m$ with $p\ge q>m$.  Let $\eta^i_1$ and $\eta^i_2$ be
odd variables for $i=1,\ldots,p$.  Then $${\rm det}
\begin{pmatrix}
x_1^1&x_1^2&\cdots&x_1^p\\
x_2^1&x_2^2&\cdots&x_2^p\\
\vdots&\vdots&\cdots&\vdots\\
x_m^1&x_m^2&\cdots&x_m^p\\
\eta_1^1&\eta_1^2&\cdots&\eta_1^p\\
\eta_1^1&\eta_1^2&\cdots&\eta_1^p\\
\vdots&\vdots&\cdots&\vdots\\
\vdots&\vdots&\cdots&\vdots\\
\eta_1^1&\eta_1^2&\cdots&\eta_1^p\\
\end{pmatrix}
\cdot{\rm det}
\begin{pmatrix}
x_1^1&x_1^2&\cdots&x_1^q\\
x_2^1&x_2^2&\cdots&x_2^q\\
\vdots&\vdots&\cdots&\vdots\\
x_m^1&x_m^2&\cdots&x_m^q\\
\eta_1^1&\eta_1^2&\cdots&\eta_1^q\\
\eta_2^1&\eta_2^2&\cdots&\eta_2^q\\
\vdots&\vdots&\cdots&\vdots\\
\eta_2^1&\eta_2^2&\cdots&\eta_2^q\\
\end{pmatrix}=0.
$$
\end{cor}
\begin{proof}
Consider $\Delta_{1,p}$ and $\Delta_{2,q}$ where $p\ge q$.  By
\thmref{glpmn-duality} the product $\Delta_{1,p}\Delta_{2,q}$ is a
highest weight vector and thus is annihilated by all operators in
\eqnref{radicalpmn}.  In particular applying the operator
$\sum_{j=1}^p\eta_{1}^{j}\frac{\partial}{\partial\eta_2^j}$ to
$\Delta_{1,p}\Delta_{2,q}$ and dividing by $(q-m)$ we obtain the
desired identity.
\end{proof}

\begin{rem}  \rm
The above corollary gives rise to identities involving minors in
even variables $x$'s by looking at the coefficient of a fixed
Grassmann monomial involving $\eta$'s. We do not know of other
direct proof of these identities.
\end{rem}

It is well known (cf.~\cite{OV}) that as a $gl(p)$-module
$S^i(\C^p)\otimes\Lambda^j(\C^p)$, for $j=1,\ldots,p$, decomposes
into a direct sum two irreducible components of highest weights
$i\epsilon_1+\sum_{k=1}^j\epsilon_k$ and
$i\epsilon_1+\sum_{k=2}^{j+1}\epsilon_k$, respectively.  We can
also get this result from \thmref{glpmn-duality} and in addition
obtain explicit formulas of the highest weight vectors. To do so
consider $S^{k}(\C^p\otimes\C^{1|1})\cong S^k(\C^{p|p})
\cong\sum_{i+j=k}S^i(\C^{p|0})\otimes\Lambda^j(\C^{0|p})$. Now
according to \thmref{glpmn-duality} all the  $gl(p)\times gl(1|1)$
highest weight vectors inside $S^{k}(\C^p\otimes\C^{1|1})$ are
given by $(x_1^1)^i\Delta_{1,j}$, where $j=1,\ldots,p$ and
$i+j=k$.  These vectors are of course $gl(p)$ highest weight
vectors.  Now a simple calculation shows that applying the
negative root vector of $gl(1|1)$ to $(x_1^1)^i\Delta_{1,j}$ we
obtain a non-zero multiple of
$(x_1^1)^i\eta_1^1\eta_1^2\ldots\eta_1^j$, while applying the
negative root vector again gives zero.  Thus the vectors
$(x_1^1)^i\Delta_{1,j}$ and
$(x_1^1)^i\eta_1^1\eta_1^2\ldots\eta_1^j$ exhaust all $gl(p)$
highest weight vectors inside the space
$\sum_{i+j=k}S^i(\C^p)\otimes\Lambda^j(\C^p)$. To conclude the
proof we observe that the vectors $(x_1^1)^{i-1}\Delta_{1,j+1}$
and $(x_1^1)^i\eta_1^1\eta_1^2\ldots\eta_1^j$ lie in
$S^i(\C^p)\otimes\Lambda^j(\C^p)$, with weights
$i\epsilon_1+\sum_{k=1}^j\epsilon_k$ and
$i\epsilon_1+\sum_{k=2}^{j+1}\epsilon_k$, respectively.

\subsection{Highest Weight Vectors: the Case $p=m$} \label{sect_pmequal}

In this section we shall find $gl(m|q)\times gl(m|n)$ highest
weight vectors that appear in the decomposition of
$\C[\x,\xibf,\etabf,\y]$. By \thmref{supergl-duality} we need to
construct a vector in $\C[\x, \xibf, \etabf, \y]$ annihilated by
all operators in \eqnref{radicalpq} and \eqnref{radicalmn} of
weight corresponding to the Young diagram $\la$
\begin{equation}\label{youngp=m}
{\beginpicture
\setcoordinatesystem units <1.5pc,1.5pc>
point at 0 3
\setplotarea x from 0 to 1.5, y from 0 to 6
\plot 0 0 1 0 1 2 2 2 2 3 3 3 3 4 4 4 4 5 7 5 7 6 0 6 0 0 /
\plot 0 3 2 3 /
\put{$\vdots$} <0pt,2pt> at 0.5 4.5
\put{$\cdots$} <0pt,2pt> at 3.5 5.3
\put{$\cdots$} <0pt,2pt> at 3 4.3
\put{$\cdots$} <0pt,2pt> at 2 3.3
\put{$\vdots$} <0pt,2pt> at 0.5 1.5
\put{$m\left\{\hbox to 0pt{
\vrule width 0 pc height 1.8pc depth 1.8pc\hfil}
  \right.$} at -.6 4.5
\put{$s\left\{\hbox to 0pt{
\vrule width 0 pc height 1.8pc depth 1.8pc\hfil}
  \right.$} at -.5 1.5
\put{$\underbrace{\hbox to 7.3pc{\ \hfill}}_{t-r}$} at 4.7 -0.5
\put{$\underbrace{\hbox to 3pc{\ \hfill}}_{r}$} at 1 -0.5
\endpicture}
\end{equation}
where $r\le min(q,n)$ (which we will always assume for this
section).

First we remark that if $\la$ has length less than or equal to $m$
then it is easy to check that a formula for the corresponding
highest weight vector is given by
$\Delta_{\la_1'}\cdots\Delta_{\la_t'}$.  So we may assume that the
length of $\la$ exceeds $m$.

As before we cut up this Young diagram into two diagrams, namely
\begin{equation}\label{splitp=m}
{\beginpicture
\setcoordinatesystem units <1.5pc,1.5pc>
point at 0 3
\setplotarea x from 0 to 1.5, y from 0 to 6
\plot 0 3 2 3 2 6 /
\plot 0 0 0 6 2 6 2 2 1 2 1 0 0 0 /
\put{$\vdots$} <0pt,2pt> at 0.5 4.5
\put{$\vdots$} <0pt,2pt> at 0.5 1.5
\put{$m\left\{\hbox to 0pt{
\vrule width 0 pc height 1.8pc depth 1.8pc\hfil}
  \right.$} at -.6 4.5
\put{$s\left\{\hbox to 0pt{
\vrule width 0 pc height 1.8pc depth 1.8pc\hfil}
  \right.$} at -.5 1.5
\put{$\underbrace{\hbox to 3pc{\ \hfill}}_{r}$} at 1 -0.5
\plot 6 3 7 3 7 4 8 4 8 5 11 5 11 6 6 6 6 3 /
\plot 0 3 2 3 /
\put{$\cdots$} <0pt,2pt> at 7.5 5.3
\put{$\cdots$} <0pt,2pt> at 7 4.3
\put{$m\left\{\hbox to 0pt{
\vrule width 0 pc height 1.8pc depth 1.8pc\hfil}
  \right.$} at 5.4 4.5
\put{$\underbrace{\hbox to 7.3pc{\ \hfill}}_{t-r}$} at 8.7 -0.5
\endpicture}
\end{equation}

Denoting the second diagram by $\mu$ and $v$ a highest weight
vector associated to the first diagram, it is easy to see that the
product $v\Delta_{\mu_1'}\cdots\Delta_{\mu_{t-r}'}$is a highest
weight vector for the diagram $\la$. Thus our task reduces to
finding a highest weight vector associated to the first diagram in
\eqnref{splitp=m}.

We claim that a highest weight vector associated to the first
diagram in \eqnref{splitp=m} can be essentially obtained by taking
a product of those associated to $s$ diagrams of rectangular shape
\begin{equation}\label{depth=m+1}
{\beginpicture
\setcoordinatesystem units <1.5pc,1.5pc>
point at 0 2
\setplotarea x from 0 to 1.5, y from 0 to 4
\plot 0 0 3 0 3 4 0 4 0 0 /
\plot 0 1 3 1 /
\put{$m\left\{\hbox to 0pt{
\vrule width 0 pc height 1.8pc depth 1.8pc\hfil}
  \right.$} at -.6 2.5
\put{$1\left\{\hbox to 0pt{
\vrule width 0 pc height 0.3pc depth 0.3pc\hfil}
  \right.$} at -.6 0.5
\put{$\underbrace{\hbox to 4.5pc{\ \hfill}}_{\la_{m+i}}$} at 1.5 -0.5
\endpicture}
\end{equation}
and dividing by a suitable power of $\Delta_{m}$.  Indeed taking
the product of two highest weight vectors for the Young diagram of
shape \eqnref{depth=m+1} of widths $\la_{m+i}$ and $\la_{m+i+1}$
respectively produces a highest weight vector for the Young
diagram
\begin{equation*}
{\beginpicture
\setcoordinatesystem units <1.5pc,1.5pc>
point at 0 2.5
\setplotarea x from 0 to 1.5, y from 0 to 5
\plot 0 0 2 0 2 1 3 1 3 5 0 5 0 0 /
\plot 0 1 2 1 /
\plot 0 2 3 2 /
\plot 3 2 5 2 5 5 3 5 /
\put{$m\left\{\hbox to 0pt{
\vrule width 0 pc height 1.8pc depth 1.8pc\hfil}
  \right.$} at -.6 3.5
\put{$2\left\{\hbox to 0pt{
\vrule width 0 pc height 0.8pc depth 0.8pc\hfil}
  \right.$} at -.6 1
\put{$\underbrace{\hbox to 3.0pc{\ \hfill}}_{\la_{m+i+1}}$} at 1
-0.6
\put{$\underbrace{\hbox to 3.0pc{\ \hfill}}_{\la_{m+i+1}}$} at 4
-0.6
\put{$\underbrace{\hbox to 4pc{\ \hfill}}_{\la_{m+i}}$} at 1.5
-1.5
\endpicture}
\end{equation*}
Once we verify that the product is non-zero, we may divide it by
$(\Delta_m)^{\la_{m+i+1}}$ and the resulting vector is a highest
weight vector for the diagram
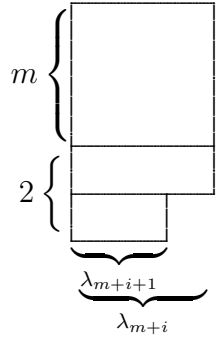
\begin{equation*}
{\beginpicture
\setcoordinatesystem units <1.5pc,1.5pc>
point at 0 2.5
\setplotarea x from 0 to 1.5, y from 0 to 5
\plot 0 0 2 0 2 1 3 1 3 5 0 5 0 0 /
\plot 0 1 2 1 /
\plot 0 2 3 2 /
\put{$m\left\{\hbox to 0pt{
\vrule width 0 pc height 1.8pc depth 1.8pc\hfil}
  \right.$} at -.6 3.5
\put{$2\left\{\hbox to 0pt{
\vrule width 0 pc height 0.8pc depth 0.8pc\hfil}
  \right.$} at -.6 1
\put{$\underbrace{\hbox to 3.0pc{\ \hfill}}_{\la_{m+i+1}}$} at 1
-0.6
\put{$\underbrace{\hbox to 4pc{\ \hfill}}_{\la_{m+i}}$} at 1.5
-1.5
\endpicture}
\end{equation*}
Similarly by taking a product of $s$ such vectors associated to
the $s$ diagrams of the form \eqnref{depth=m+1} of widths
$\la_{m+1},\ldots, \la_{m+s}$, respectively, and dividing by
$\Delta_m^{\la_{m+2}+\ldots+\la_{m+s}}$ we obtain a highest weight
vector associated to the first diagram of \eqnref{splitp=m}. So
our task now is to find a formula for a highest weight vector
corresponding to a Young diagram of shape \eqnref{depth=m+1}.
(From the explicit formula it will follow immediately that a
product of $s$ vectors of such type is non-zero).

Let us put $r=\la_{m+i}$ in \eqnref{depth=m+1}. We define the
$r\times r$ matrix $Y$ and the $m\times m$ matrix $X$ as follows:
\begin{equation}\label{xy}
Y:=
\begin{pmatrix}
y_1^1&y_1^2&\cdots&y_1^r\\
y_2^1&y_2^2&\cdots&y_2^r\\
\vdots&\vdots&\cdots&\vdots\\
y_r^1&y_r^2&\cdots&y_r^r\\
\end{pmatrix},\quad
X:=
\begin{pmatrix}
x_1^1&x_1^2&\cdots&x_1^m\\
x_2^1&x_2^2&\cdots&x_2^m\\
\vdots&\vdots&\cdots&\vdots\\
x_m^1&x_m^2&\cdots&x_m^m\\
\end{pmatrix}.
\end{equation}

Given a Young diagram $\la$ of rectangular shape (see
\eqnref{depth=m+1}) consisting of $m$ rows and $r$ columns, we
consider {\em marked diagrams} $D$ obtained by marking the boxes
in $\la$ subject to the restriction that each column can contain
no more than one marked box. For example the following is a marked
diagram in the case $r=6$ and $m=4$:
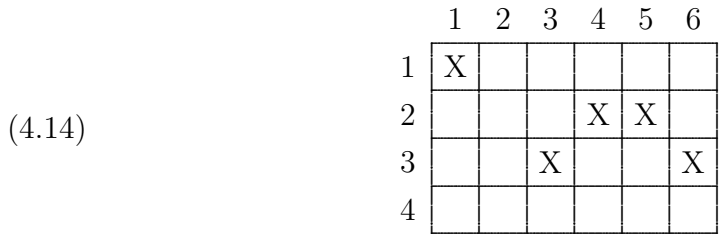
\begin{equation}\label{example}
{\beginpicture
\setcoordinatesystem units <1.5pc,1.5pc>
point at 0 2
\setplotarea x from 0 to 1.5, y from 0 to 4
\plot 0 0 6 0 6 4 0 4 0 0 /
\plot 0 1 6 1 /
\plot 0 2 6 2 /
\plot 0 3 6 3 /
\plot 1 0 1 4 /
\plot 2 0 2 4 /
\plot 3 0 3 4 /
\plot 4 0 4 4 /
\plot 5 0 5 4 /
\put {$1$} at 0.5 4.5
\put {$2$} at 1.5 4.5
\put {$3$} at 2.5 4.5
\put {$4$} at 3.5 4.5
\put {$5$} at 4.5 4.5
\put {$6$} at 5.5 4.5
\put {$1$} at -0.5 3.5
\put {$2$} at -0.5 2.5
\put {$3$} at -0.5 1.5
\put {$4$} at -0.5 0.5
\put {X} at 0.5 3.5
\put {X} at 2.5 1.5
\put {X} at 4.5 2.5
\put {X} at 5.5 1.5
\put {X} at 3.5 2.5
\endpicture}
\end{equation}
To each such a marked diagram $D$ we may associate an $r\times r$
matrix $Y_D$ obtained from $Y$ as follows. For each marked box,
say in the $i$-th column and $j$-th row, we replace the $i$-th row
of the matrix $Y$ by the vector
$(\xi_j^1,\xi_j^2,\ldots,\xi_j^r)$. The resulting matrix will be
denoted by $Y_D$.  For instance in our example \eqnref{example}
the matrix $Y_D$ is
\begin{equation*}
Y_D =
\begin{pmatrix}
\xi_1^1&\xi_1^2&\xi_1^3&\xi_1^4&\xi_1^5&\xi_1^6\\
y_1^1&y_1^2&y_1^3&y_1^4&y_1^5&y_1^6\\
\xi_3^1&\xi_3^2&\xi_3^3&\xi_3^4&\xi_3^5&\xi_3^6\\
\xi_2^1&\xi_2^2&\xi_2^3&\xi_2^4&\xi_2^5&\xi_2^6\\
\xi_2^1&\xi_2^2&\xi_2^3&\xi_2^4&\xi_2^5&\xi_2^6\\
\xi_3^1&\xi_3^2&\xi_3^3&\xi_3^4&\xi_3^5&\xi_3^6\\
\end{pmatrix}.
\end{equation*}
To each such diagram $D$ we may also associate $r$ $m\times m$
matrices $X_i$ $(i=1,\ldots,r)$ obtained from $X$ as follows. If
the $i$-th column of $D$ is not marked, then $X_i=X$.  If the
$i$-th column is marked at the $j$-th row, then $X_i$ is the
matrix obtained from $X$ by replacing its $j$-th row by the vector
$(\eta_i^1,\eta_i^2,\ldots,\eta_i^m)$.  As an illustration, the
diagram in our example \eqnref{example} gives rise to the matrices
\begin{equation*}
X_1=
\begin{pmatrix}
\eta_1^1&\eta_1^2&\eta_1^3&\eta_1^4\\
x_2^1&x_2^2&x_2^3&x_2^4\\
x_3^1&x_3^2&x_3^3&x_3^4\\
x_4^1&x_4^2&x_4^3&x_4^4\\
\end{pmatrix},
X_2=X,
X_3=\begin{pmatrix}
x_1^1&x_1^2&x_1^3&x_1^4\\
x_2^1&x_2^2&x_2^3&x_2^4\\
\eta_3^1&\eta_3^2&\eta_3^3&\eta_3^4\\
x_4^1&x_4^2&x_4^3&x_4^4\\
\end{pmatrix}\ {\rm etc.}
\end{equation*}

Let $|D|$ denote the total number of marked boxes in the diagram
$D$. Set $$\Delta_D:={\rm det}X_D{\rm det}Y_D, $$ where by ${\rm
det}X_D$ we mean $\prod_{i=1}^r{\rm det}X_i$ arranged in
increasing order.  We can now state the following theorem.

\begin{thm}\label{mainthmp=m}
The vector $\sum_{D}(-1)^{\frac{1}{2}|D|(|D|-1)}\Delta_D$ is a
$gl(m|q)\times gl(m|n)$ highest weight vector in $\C [\x, \xibf,
\etabf, \y ]$ corresponding to the rectangular Young diagram of
length $m +1$ and width $r$, where the summation over $D$ ranges
over all possible marked $m \times r$ diagrams.
\end{thm}

\begin{proof}
We first show that $\sum_{D}(-1)^{\frac{1}{2}|D|(|D|-1)}\Delta_D$ indeed has
the correct weight.

First note that diagram \eqnref{depth=m+1} corresponds to the
$gl(m|q)\times gl(m|n)$-weight $\sum_{i=1}^m
r\epsilon_i+\sum_{i=1}^m
r\tilde{\epsilon}_i+\sum_{j=1}^r\delta_j+\sum_{j=1}^r\tilde{\delta}_j$.
Let $\tilde{D}_i$, $j=1,\ldots,m$, denote the $m$ disjoint subsets
of $\{1,\ldots,r\}$ defined by the condition that
$j\in\tilde{D}_i$ if and only if $D$ contains a marked box at its
$j$-th column and $i$-th row. Put
$\tilde{D}=\cup_{i=1}^m\tilde{D}_i$ and
$\tilde{D}^c=\{1,\ldots,r\}-\tilde{D}$.  The weight of ${\rm
det}Y_D$ is $\sum_{j\in\tilde{D}^c}\delta_j
+\sum_{i=1}^m|\tilde{D}_i|\epsilon_i+\sum_{j=1}^r\tilde{\delta}_j$.
Now ${\rm det}X_D$ has weight
$r\sum_{i=1}^m\epsilon_i-\sum_{i=1}^m|\tilde{D}_i|\epsilon_i+
\sum_{j\in\tilde{D}}\delta_j+r\sum_{i=1}^m\tilde{\epsilon}_i$.
Hence each ${\rm det}Y_D{\rm det}X_D$ has weight $\sum_{i=1}^m
r\epsilon_i+\sum_{i=1}^m
r\tilde{\epsilon}_i+\sum_{j=1}^r\delta_j+\sum_{j=1}^r\tilde{\delta}_j$,
as required.

Hence by Corollary~\ref{cor_hwt} it is sufficient to show that
\eqnref{radicalmn} annihilates it, namely
\begin{eqnarray}
\left( \sum_{j=1}^m x_{s-1}^{j}\frac{\partial}{\partial
x_{s}^j}+\sum_{j=1}^q \xi_{s-1}^{j}\frac{\partial}{\partial
\xi_{s}^j}\right)
(\sum_{D}(-1)^{\frac{1}{2}|D|(|D|-1)}\Delta_D ) & =& 0,
\label{step1}\\
\left(
\sum_{j=1}^m\eta_{s-1}^{j}\frac{\partial}{\partial\eta_s^j}+
\sum_{j=1}^q y_{s-1}^{j}\frac{\partial}{\partial y_s^j} \right)
(\sum_{D}(-1)^{\frac{1}{2}|D|(|D|-1)}\Delta_D ) & =& 0 ,
\label{step2}\\
\left( \sum_{j=1}^m
x_m^j\frac{\partial}{\partial\eta_1^j} - \sum_{j=1}^q
\xi_m^j\frac{\partial}{\partial y_1^j} \right)
(\sum_{D}(-1)^{\frac{1}{2}|D|(|D|-1)}\Delta_D ) & =& 0.\label{step3}
\end{eqnarray}

We will first establish \eqnref{step1}. Note that the simple root
vector $\sum_{j=1}^m x_{s-1}^{j}\frac{\partial}{\partial
x_{s}^j}+\sum_{j=1}^q \xi_{s-1}^{j}\frac{\partial}{\partial
\xi_{s}^j}$ maps the vectors $(x_s^1,\ldots,x_s^m)$ to
$(x_{s-1}^1,\ldots,x_{s-1}^m)$ and $(\xi_s^1,\ldots,\xi_s^q)$ to
$(\xi_{s-1}^1,\ldots,\xi_{s-1}^q)$.  For a diagram $D$, let us
denote by $D_{\uparrow_{s}^{s-1}}$ a diagram obtained from $D$ by
moving each marked box in its $s$-th row to the box above it in
the $s-1$-st row.  Analogously we define
$D_{\downarrow_{s}^{s-1}}$ a diagram obtained from $D$ by moving
each marked box in the $s-1$-st row to the box below it in the
$s$-th row.  It is easy to check
\begin{align*}
\left( \sum_{j=1}^m x_{s-1}^{j}\frac{\partial} {\partial x_{s}^j}
+\sum_{j=1}^q \xi_{s-1}^{j} \frac{\partial}{\partial \xi_{s}^j}
\right) & (\Delta_D) = \\
 & \sum_{D_{\uparrow^{s-1}_s}}{\rm det}X_D{\rm det}
Y_{D_{\uparrow^{s-1}_s}} - \sum_{D_{\downarrow_{s}^{s-1}}} {\rm
det} X_{D_{\downarrow_{s}^{s-1}}} {\rm det}Y_D.
\end{align*}
Thus we have
\begin{align}
\left( \sum_{j=1}^m x_{s-1}^{j}\frac{\partial}{\partial x_{s}^j}
\right. &\left. +\sum_{j=1}^q \xi_{s-1}^{j}
\frac{\partial}{\partial \xi_{s}^j} \right)
(\sum_{|D|=k}\Delta_D) =  \nonumber\\
 &\sum_{|D|=k}\left(
\sum_{D_{\uparrow^{s-1}_s}}{\rm det}X_D{\rm
det}Y_{D_{\uparrow^{s-1}_s}}- \sum_{D_{\downarrow_{s}^{s-1}}}{\rm
det}X_{D_{\downarrow_{s}^{s-1}}}{\rm det}Y_D \right). \label{aux1}
\end{align}
But evidently $\sum_{|D|=k} {\rm det}X_{D_{\downarrow^{s-1}_s}}
{\rm det}Y_D = \sum_{|D|=k} {\rm det}X_D{\rm det}
Y_{D_{\uparrow^{s-1}_s}}$ thanks to the equality
$(D_{\uparrow^{s-1}_s})_{\downarrow^{s-1}_s} =D$. Hence the
right-hand side of \eqnref{aux1} is zero, proving \eqnref{step1}.

Our next step is to prove \eqnref{step2}.  In this case
$\sum_{j=1}^m\eta_{s-1}^{j}\frac{\partial}{\partial\eta_s^j}+
\sum_{j=1}^q y_{s-1}^{j}\frac{\partial}{\partial y_s^j}$ maps the
vectors $(\eta_s^1,\ldots,\eta_s^m)$ to
$(\eta_{s-1}^1,\ldots,\eta_{s-1}^m)$ and $(y_s^1,\ldots,y_s^q)$ to
$(y_{s-1}^1,\ldots,y_{s-1}^q)$.  For a diagram $D$ such that $j\in
\tilde{D}_i$ and $l\not\in\tilde{D}_i$ we denote by
$D_{j\rightarrow l}$ the diagram obtained from $D$ by removing $j$
from $\tilde{D}_i$ and adding $l$ to $\tilde{D}_i$. For a fixed
$k$ we write
\begin{equation*}
\sum_{|D|=k}\Delta_D =\sum_{|D|=k}\left( \sum_{s,s-1\in
\tilde{D}}\Delta_D + \sum_{s\not\in\tilde{D}{\rm\ or\
}s-1\not\in\tilde{D}}\Delta_D \right).
\end{equation*}

First observe that
\begin{align*}
&\left(\sum_{j=1}^m\eta_{s-1}^{j}\frac{\partial}{\partial\eta_s^j}+
\sum_{j=1}^q y_{s-1}^{j}\frac{\partial}{\partial y_s^j} \right)
(\sum_{|D|=k}\sum_{s,s-1\in \tilde{D}}\Delta_D ) \\ =
&\sum_{|D|=k}\sum_{s,s-1\in \tilde{D}} \left( \left(
\sum_{j=1}^m\eta_{s-1}^{j}\frac{\partial}{\partial\eta_s^j}+
\sum_{j=1}^q y_{s-1}^{j}\frac{\partial}{\partial y_s^j} \right)
({\rm det}X_D) \right){\rm det}Y_D \\ =& \ 0.
\end{align*}
This is because if $s,s-1\in\tilde{D}_i$, for some $i$, then the
term
\begin{align*}
& \left( \sum_{j=1}^m
\eta_{s-1}^{j}\frac{\partial}{\partial\eta_s^j} + \sum_{j=1}^q
y_{s-1}^{j} \frac{\partial}{\partial y_s^j}
   \right) ({\rm det}X_D)\\
  =& \cdots {\rm det}X_{s-1}(i){\rm det}X_{s-1}(i)\cdots  \\
  =& \ 0,
\end{align*}
where in general $X_b(a)$ is the matrix obtained from $X$ by
replacing the $a$-th row with the vector $(\eta_{b}^1, \ldots,
\eta_{b}^m)$.

Now if $s\in \tilde{D}_i$ and $s-1\in\tilde{D}_l$, then
$$ \left(
\sum_{j=1}^m\eta_{s-1}^{j}\frac{\partial}{\partial\eta_s^j} +
\sum_{j=1}^q y_{s-1}^{j}\frac{\partial}{\partial y_s^j} \right)
({\rm det}X_D) =\cdots{\rm det}X_{s-1}(l){\rm
det}X_{s-1}(i)\cdots. $$
Let $D'$ be the same diagram as $D$, except $s\in \tilde{D'}_l$
and $s-1\in\tilde{D'}_i$.  Then we have
$$
\left(\sum_{j=1}^m\eta_{s-1}^{j}\frac{\partial}{\partial\eta_s^j}+
\sum_{j=1}^q y_{s-1}^{j}\frac{\partial}{\partial y_s^j} \right)
({\rm det}X_{D'}) =\cdots{\rm det}X_{s-1}(i){\rm
det}X_{s-1}(l)\cdots. $$
Of course $Y_D=Y_{D'}$ and ${\rm det}X_{s-1}(i)$ anticommutes with
$ {\rm det}X_{s-1}(l)$, so
$$ \left(
\sum_{j=1}^m\eta_{s-1}^{j}\frac{\partial}{\partial\eta_s^j} +
\sum_{j=1}^q y_{s-1}^{j}\frac{\partial}{\partial y_s^j}
 \right)
({\rm det}X_D{\rm det}Y_D+{\rm det}X_{D'}{\rm det}Y_{D'}) =0.$$
Next we observe that if $D$ is a diagram such that
$s,s-1\not\in\tilde{D}$ then
\[
\left( \sum_{j=1}^m\eta_{s-1}^{j}\frac{\partial}{\partial\eta_s^j}
+ \sum_{j=1}^q y_{s-1}^{j} \frac{\partial}{\partial y_s^j} \right)
({\rm det}X_D {\rm det}Y_D) =0 ,
\]
so that our task of proving \eqnref{step2} reduces to proving that

\begin{align} \label{aux2}
\left( \sum_{j=1}^m\eta_{s-1}^{j}\frac{\partial}{\partial\eta_s^j}
+ \sum_{j=1}^q y_{s-1}^{j}\frac{\partial}{\partial y_s^j} \right)
 & \left( \sum_{\stackrel{s\in\tilde{D}
}{s-1\not\in\tilde{D} } } {\rm det}X_D{\rm det}Y_D
+\sum_{\stackrel{s-1\in\tilde{D} }{ s\not\in\tilde{D}}} {\rm
det}X_D{\rm det}Y_D \right)   \nonumber \\ =& 0,
\end{align}
where the sum is over all diagrams $D$ with $|D|=k$.  But the
left-hand side of \eqnref{aux2} is equal to
\begin{equation*}
\sum_{s\in\tilde{D},s-1\not\in\tilde{D}}{\rm
det}X_{D_{s\rightarrow s-1}}{\rm
det}Y_D-\sum_{s-1\in\tilde{D},s\not\in\tilde{D}}{\rm det}X_D{\rm
det}Y_{D_{s-1\rightarrow s}} =0.
\end{equation*}

To complete the proof we now need to verify \eqnref{step3}.  The
odd simple root vector $\sum_{j=1}^m
x_m^j\frac{\partial}{\partial\eta_1^j}- \sum_{j=1}^q
\xi_m^j\frac{\partial}{\partial y_1^j}$ has the effect of changing
the vectors $(\eta_1^1,\ldots,\eta_1^m)$ to $(x_m^1,\ldots,x_m^m)$
and $(y_1^1,\ldots,y_1^q)$ to $-(\xi_m^1,\ldots,\xi_m^q)$.  If $D$
is a diagram such that $1\in D_j$ with $j\not=m$, then
\[
\left( \sum_{j=1}^m x_m^j\frac{\partial}{\partial\eta_1^j} -
\sum_{j=1}^q \xi_m^j\frac{\partial}{\partial y_1^j} \right) \cdot
 \Delta_D =0.
\]
Thus
\begin{align}
&\left( \sum_{j=1}^m  x_m^j\frac{\partial}{\partial\eta_1^j}-
\sum_{j=1}^q \xi_m^j \frac{\partial}{\partial y_1^j}\right)
(\sum_{D}(-1)^{\frac{1}{2}|D|(|D|-1)}\Delta_D) = \nonumber \\
 & \left( \sum_{j=1}^m x_m^j\frac{\partial}{\partial\eta_1^j} -
\sum_{j=1}^q \xi_m^j\frac{\partial}{\partial y_1^j} \right)
(\sum_{D,1\in\tilde{D}_m}(-1)^{\frac{1}{2}|D|(|D|-1)}\Delta_D +
\sum_{D,1\not\in\tilde{D}}(-1)^{\frac{1}{2}|D|(|D|-1)}
\Delta_D).\label{aux3}
\end{align}

For a diagram $D$ with $1\not\in\tilde{D}$ (resp. with
$1\in\tilde{D}_m$) we denote by $D^{+}$ (resp. $D^-$) the diagram
obtained from $D$ by adding $1$ to $\tilde{D}_m$ (resp. by
removing $1$ from $\tilde{D}$). Then \eqnref{aux3} becomes
\begin{equation}\label{aux4}
\sum_{D,1\in\tilde{D}_m}(-1)^{{\frac{1}{2}}|D|(|D|-1)}{\rm det}X_{D^-}{\rm
det}Y_D-
\sum_{D,1\not\in\tilde{D}}(-1)^{\frac{1}{2}|D|(|D|-1)+|D|}{\rm det}X_D{\rm
det}Y_{D^+}.
\end{equation}
Setting $D'=D^-$ in the first sum of \eqnref{aux4} we may rewrite
\eqnref{aux4} as
\begin{align*}
\sum_{D',1\not\in\tilde{D'}}(-1)^{{\frac{1}{2}}|D'|(|D'|+1)} &{\rm
det}X_{D'} \ {\rm det}Y_{{D'}^+} \\ - &
\sum_{D,1\not\in\tilde{D}}(-1)^{\frac{1}{2}|D|(|D|-1)+|D|}{\rm
det}X_D\ {\rm det}Y_{D^+}=0.
\end{align*}
\end{proof}

We will denote the vector
$\sum_{D}(-1)^{\frac{1}{2}|D|(|D|-1)}\Delta_D$ by $\Gamma_r$. It
is clear that a product of $\Gamma_r$'s (not necessary for the
same value $r$) remains nonzero.  Thus a highest weight vector for
an arbitrary Young diagram of shape \eqnref{youngp=m} can be
constructed using such vectors, as described earlier in this
section. We summarize the results in this section in the following
theorem.

\begin{thm} An irreducible representation $V^{\la}_{m|q}\otimes V^{\la}_{m|n}$
of $gl(m|q)\times gl(m|n)$  appears in the decomposition of
$\C[\x, \xibf, \etabf, \y]$ if and only if $\la$ is associated to
a Young diagram with $\la_{m+1}\le {\rm min}(q,n)$. Let $t$ be the
length of $\la'$. Then
\begin{itemize}
\item[(1)] if the length of $\la$ does not exceed $m$, then a highest weight
vector is given by $\Delta_{\la'_1}\cdots\Delta_{\la'_t}$.
\item[(2)] if the length of $\la$ is $m+s$, $s\ge 1$, let $0\le r\le {\rm
min}(q,n)$ be such that $\la'_{r}> m$ and $\la'_{r+1}\le m$, then
a highest weight vector corresponding to $\la$ is given by

\begin{equation}\label{formulap=m}
(\Delta_m)^{- (\la_{m+2} + \ldots+\la_{m+s}) } \Gamma_{\la_{m+1}}
\Gamma_{\la_{m+2}} \cdots \Gamma_{\la_{m+s}} \Delta_{\la'_{r+1}}
\cdots \Delta_{\la'_{t}}.
\end{equation}
\end{itemize}
\end{thm}

We will obtain a more explicit formula for \eqnref{formulap=m} in
the next section.
\subsection{Highest Weight Vectors: the General Case} \label{sect_general}

We consider now the general case. Without loss of generality we
may assume $p\ge m$.

According to \thmref{supergl-duality} an irreducible
$gl(p|q)\times gl(m|n)$-module $V^\la_{p|q}\otimes V^\la_{m|n}$
appears in the decomposition of $\C[\x, \xibf, \etabf, \y]$ if and
only if $\la_{m+1}\le n$ and $\la_{p+1}\le q$. If the length of
the Young diagram $\la$ is less than or equal to $m$, then
$\Delta_{\la'_1}\cdots\Delta_{\la'_t}$ is the desired highest
weight vector, where $t$ is the length of $\la'$. If the length of
$\la$ exceeds $m$, but is less than or equal to $p$, then we see
that the vector given in \thmref{glpmn-duality} provides a formula
for the highest weight vector in this case as well.  Thus it
remains to study the case when the length of $\la$ exceeds $p$.

So we are to consider a Young diagram of the form:
\begin{equation}\label{youngp>m}
{\beginpicture
\setcoordinatesystem units <1.5pc,1.5pc>
point at 0 3
\setplotarea x from 0 to 1.5, y from 0 to 7
\plot 0 0 1 0 1 2 2 2 2 3 3 3 3 4 4 4 4 5 5 5 5 6 7 6 7 7 0 7 0 0 /
\plot 0 5 4 5 /
\plot 0 3 2 3 /
\put{$\vdots$} <0pt,2pt> at 0.5 6
\put{$\cdots$} <0pt,2pt> at 3 6
\put{$\cdots$} <0pt,2pt> at 2 4.0
\put{$\vdots$} <0pt,2pt> at 0.5 2.0
\put{$m\left\{\hbox to 0pt{
\vrule width 0 pc height 1.3pc depth 1.3pc\hfil}
  \right.$} at -.6 6
\put{$p-m\left\{\hbox to 0pt{
\vrule width 0 pc height 1.0pc depth 1.0pc\hfil}
  \right.$} at -1.2 4
\put{$s\left\{\hbox to 0pt{
\vrule width 0 pc height 1.5pc depth 1.5pc\hfil}
  \right.$} at -.6 1.5
\put{$\underbrace{\hbox to 3pc{\ \hfill}}_{r'-r}$} at 3 -0.5
\put{$\underbrace{\hbox to 5pc{\ \hfill}}_{t-r'}$} at 6 -0.5
\put{$\underbrace{\hbox to 3pc{\ \hfill}}_{r}$} at 1 -0.5
\endpicture}
\end{equation}

In the case when $q\le n$, the numbers $r, r'$ satisfying the
conditions $0\le r\le q$, $0\le r'\le n$ and $r\le r'$ are
determined as follows: $\la'_{r}>p$ and $\la'_{r+1}\le p$,
$\la'_{r'}>m$ and $\la'_{r'+1}\le m$.  In the case when $q\ge n$,
the numbers $r, r'$ satisfying $0\le r\le r'\le n$ are defined in
exactly the same way. In either case we may split
\eqnref{youngp>m} into three diagrams:
\begin{equation}\label{splityoungn>q}
{\beginpicture
\setcoordinatesystem units <1.5pc,1.5pc>
point at 0 3
\setplotarea x from 0 to 15, y from 0 to 7
\plot 0 0 1 0 1 2 2 2 2 7 0 7 0 0 /
\plot 5 3 6 3 6 4 7 4 7 7 5 7 5 3 /
\plot 10 5 11 5 11 6 13 6 13 7 10 7 10 5 /
\plot 0 5 2 5 /
\plot 0 3 2 3 /
\plot 5 5 7 5 /
\put{$m\left\{\hbox to 0pt{
\vrule width 0 pc height 1.3pc depth 1.3pc\hfil}
  \right.$} at -.6 6
\put{$p-m\left\{\hbox to 0pt{
\vrule width 0 pc height 1.0pc depth 1.0pc\hfil}
  \right.$} at -1.2 4
\put{$m\left\{\hbox to 0pt{
\vrule width 0 pc height 1.3pc depth 1.3pc\hfil}
  \right.$} at 4.3 6
\put{$p-m\left\{\hbox to 0pt{
\vrule width 0 pc height 1.0pc depth 1.0pc\hfil}
  \right.$} at 3.8 4
\put{$s\left\{\hbox to 0pt{
\vrule width 0 pc height 1.5pc depth 1.5pc\hfil}
  \right.$} at -.6 1.5
\put{$m\left\{\hbox to 0pt{
\vrule width 0 pc height 1.3pc depth 1.3pc\hfil}
  \right.$} at 9.3 6
\put{$\underbrace{\hbox to 3pc{\ \hfill}}_{r'-r}$} at 6 -0.5
\put{$\underbrace{\hbox to 5pc{\ \hfill}}_{t-r'}$} at 12 -0.5
\put{$\underbrace{\hbox to 3pc{\ \hfill}}_{r}$} at 1 -0.5
\endpicture}
\end{equation}
We associate the vectors
$\Delta_{r+1,\la'_{r+1}}\cdots\Delta_{r',\la'_{r'}}$ and
$\Delta_{\la'_{r'+1}}\cdots\Delta_{\la'_t}$ to the second and
third diagrams in \eqnref{splityoungn>q} respectively. Below we
will construct a highest weight vector for the first diagram in
\eqnref{splityoungn>q}.  From the formula it will be easy to see
that the product of these three vectors is a highest weight vector
for the Young diagram \eqnref{youngp>m}.

The above discussion thus reduces our question to
finding a $gl(p|q)\times gl(m|n)$-highest weight vector
corresponding to a Young diagram of type:
\begin{equation}\label{mainp>m}
{\beginpicture
\setcoordinatesystem units <1.5pc,1.5pc>
point at 0 3
\setplotarea x from 0 to 1.5, y from 0 to 7
\plot 0 0 1 0 1 2 2 2 2 7 0 7 0 0 /
\plot 0 5 2 5 /
\plot 0 3 2 3 /
\put{$m\left\{\hbox to 0pt{
\vrule width 0 pc height 1.3pc depth 1.3pc\hfil}
  \right.$} at -.6 6
\put{$p-m\left\{\hbox to 0pt{
\vrule width 0 pc height 1.0pc depth 1.0pc\hfil}
  \right.$} at -1.2 4
\put{$s\left\{\hbox to 0pt{
\vrule width 0 pc height 1.5pc depth 1.5pc\hfil}
  \right.$} at -.6 1.5
\put{$\underbrace{\hbox to 3pc{\ \hfill}}_{r}$} at 1 -0.5
\endpicture}
\end{equation}
where $r\le{\rm min}(q,n)$. The difficulty of finding a highest
weight vector associated to such a diagram lies in the fact that
the highest weight vectors in $\C[\x, \xibf, \etabf, \y]$ no
longer form a semigroup in general.

We will now outline our strategy. We need to find highest weight
vectors in $\C[\x, \xibf, \etabf, \y]$, annihilated by
\eqnref{radicalmn} and having weight corresponding to the Young
diagram in \eqnref{mainp>m}. Recall that $\x$ denotes the set of
even variables $\{x_l^i|1\le i\le p,1\le l\le m\}$.  We introduce
a new set of even variables $x_l^i,$ $1\le i \le p,  m < l\le p$,
and denote by $\x' =\{x_l^i| 1\le i, l\le p\}$, the union of our
old set with this new set. We shall construct certain vectors in
$\C[\x', \xibf, \etabf, \y]$, which can been shown, using our
results in the previous section, that they are annihilated by
\eqnref{radicalmn}. A priori these vectors lie in $\C[\x',
\xibf, \etabf, \y]$ so that such vectors do not make sense.
However, we will show that these vectors, after dividing by a
suitable power of the determinant of the $p\times p$ matrix
$(x_l^i)$, are in fact independent of the variables $\{x_l^i|1\le
i\le p,m+1\le l\le p\}$, and thus lie in $\C[\x, \xibf, \etabf,
\y]$.

Consider a marked diagram having $m$ rows and $r$ columns with at most $r$
marked boxes
subject to the constraint that at most one marked box appears on
each column. To such a diagram $D$ we have associated in the
previous section a matrix $Y_D$, which is obtained from the
$r\times r$ matrix $Y$ (see \eqnref{xy}) by suitably replacing its
rows.  To each such diagram we now associate $p\times p$ matrices
$X_i$, for $1\le i\le r$, similar to the ones in the previous
section:  Let $X$ denote the $p\times p$ matrix:
\begin{equation}\label{Xp}
X:=\begin{pmatrix}
x_1^1&x_1^2&\cdots&x_1^p\\
x_2^1&x_2^2&\cdots&x_2^p\\
\vdots&\vdots&\cdots&\vdots\\
x_m^1&x_m^2&\cdots&x_m^p\\
x_{m+1}^1&x_{m+1}^2&\cdots&x_{m+1}^p\\
\vdots&\vdots&\cdots&\vdots\\
x_p^1&x_1^2&\cdots&x_p^p\\
\end{pmatrix}.
\end{equation}
If the $i$-th column of $D$ is not marked, then $X_i=X$.  If the
$i$-th column is marked at the $j$-th row, then $X_i$ is the
matrix obtained from $X$ by replacing its $j$-th row by the vector
$(\eta_i^1,\eta_i^2,\ldots,\eta_i^p)$.  Note that none of its
$m+1$-st to $p$-th rows are replaced. As before we define ${\rm
det}X_D:=\prod_{i=1}^r {\rm det}X_i$ (arranged in increasing
order) and define $\Gamma_r =\sum_{D}(-1)^{\frac{1}{2}|D|(D|-1)}
{\rm det}X_D {\rm det}Y_D$.

The proof of \thmref{mainthmp=m} carries over word for word to
prove

\begin{prop}
The vector $\Gamma_r$ is annihilated by \eqnref{radicalmn}.
\end{prop}

In the case when $p=m$ a highest weight vector for the first Young
diagram in \eqnref{splitp=m} is obtained essentially by taking
product of highest weight vectors for the Young diagram of type
\eqnref{depth=m+1}.  In the case when $p>m$ one can verify that
this procedure cannot be carried out, as such products are
necessarily zero.  Hence in this case we will need to find a
general formula for the Young diagram $\la$ of shape
\eqnref{mainp>m}.  To do so we will first consider matrices that
will play the same role in the case of $p\ge m$ as the $X_i$'s
play in the case $p=m$. As we will generalize diagrams to include
those that allow more than one marked box on each column, we are
led to study combinatorial identities of determinants of matrices
obtained from $X$ that have more than one row replaced by an odd
vector.  This leads us to define the following types of
determinants.

Let $X$ be the $p\times p$ matrix as in \eqnref{Xp} and let
$(\eta_j^1,\ldots,\eta_j^p)$ be an odd vector.  Let $I$ be a
subset of $\{1,\ldots,p\}$ and define $X_j(I)$ to be the matrix
obtained from $X$ by replacing its $i$-th row by the vector
$(\eta_j^1,\ldots,\eta_j^p)$, for all $i\in I$. If $I = \{ i_1,
\ldots, i_l \}$, we write $X_j(I) = X_j(i_1, \ldots, i_l)$ as
well.

\begin{lem}\label{keylemma}
We have
\begin{equation}\label{keyequation}
{\rm det}X_1(1){\rm det}X_1(2)\cdots {\rm det}X_1(p)
=\frac{1}{p!}({\rm det}X)^{p-1}{\rm det}X_1(1,\ldots,p).
\end{equation}
\end{lem}

\begin{proof}
Denoting by $R(X)$ and $L(X)$ the right-hand side and the
left-hand side of \eqnref{keyequation}, respectively, we may
regard $R$ and $L$ as functions of $X$. Since the group $GL(p)$
acts on $X$, the space of $p\times p$ matrices, by left
multiplication, it acts on functions of $X$.  To be more precise
if $A\in GL(p)$, then $(A\cdot L)(X):=L(A^{-1}X)$ and $(A\cdot
R)(X):=R(A^{-1}X)$.  We want to study the effect of this action on
$R$ and $L$.  In order to do so, consider first the action of the
three kinds of elementary matrices on them. Namely, those that
interchanges any two rows, that multiplies a row by a scalar, and
those that add a scalar multiple of a row to another.  It is
subject to a direct verification that if $A$ is any of the three
types of elementary matrices, we have
\begin{equation}\label{aux5}
(A\cdot R)(X)=({\rm det}A)^{1-p}R(X),\quad (A\cdot L)(X)=({\rm
det}A)^{1-p}L(X).
\end{equation}
Since every element in $GL(p)$ is a product of elementary
matrices, we conclude that \eqnref{aux5} holds for every $A\in
GL(p)$ as well.  Putting $X=1_p$, the identity $p\times p$ matrix, we see
that $R(1_p)=L(1_p)$ so that $R(A)=L(A)$ by \eqnref{aux5} for all
$A\in GL(p)$.  As $GL(p)$ is a Zariski open set in the space of
$p\times p$ matrices we have $R(X)=L(X)$ for any $p\times
p$ matrix $X$. \end{proof}

\begin{rem}
An alternative proof of the above lemma can be given as follows.
Let $X_{ij}$ denote the $(i,j)$-th minor of $X$. It is known that
${\rm det}(X_{ij})={\rm det}X^{p-1}$ which follows directly from a
form of the Cramer's formula $X (X_{ij}) = ({\rm det}X) 1_p$.  The
above lemma follows from this identity by expanding each
determinant in the left-hand side of \eqnref{keyequation} by the
row $(\eta_1^1,\ldots,\eta_1^p)$ and noting that ${\rm
det}X_1(1,\ldots,p)$ is equal to $ p!\  \eta_1^1 \eta_1^2 \cdots
\eta_1^p$.
\end{rem}

\begin{cor}\label{maincor}
Let $I=\{i_1,\ldots,i_{|I|}\}$ and $J=\{j_1,\ldots,j_{|J|}\}$ be
two subsets of $\{1,\ldots,p\}$ arranged in increasing order. Then
\begin{itemize}
\item[(i)] ${\rm det}X_1(I){\rm det}X_1(J)=0$ if and only if $I\cap
J\not=\emptyset$.
\item[(ii)] ${\rm det}X_1(i_1){\rm det}X_1(i_2)\cdots {\rm
det}X_1(i_{|I|})=\frac{1}{|I|!}({\rm det}X)^{|I|-1}{\rm det}X_1(I).$
\item[(iii)] For $I\cap J=\emptyset$ we have
$${\rm det}X_1(I){\rm
det}X_1(J)=\epsilon_{IJ}\frac{|I|!|J|!}{(|I|+|J|)!}{\rm det}X{\rm
det}X_1(I\cup J),$$ where $\epsilon_{IJ}$ is the sign of the
permutation that arranges the ordered tuple
$(i_1,\ldots,i_{|I|},j_1,\ldots,j_{|J|})$ in increasing order.
\end{itemize}
\end{cor}
\begin{proof}
(i) is an obvious consequence of \lemref{keylemma}.

For (ii) let $I^c=\{k_1,\ldots,k_{|I^c|}\}$ denote the
complementary subset of $I$ in $\{1,\ldots,p\}$ put in increasing
order. We apply successively the differential operators
$\sum_{j=1}^px^j_{k_1}\frac{\partial}{\partial\eta_1^j}$,
$\sum_{j=1}^px^j_{k_2}\frac{\partial}{\partial\eta_1^j}$,
$\ldots$,
$\sum_{j=1}^px^j_{k_{|I^c|}}\frac{\partial}{\partial\eta_1^j}$ to
\eqnref{keyequation} and find that $$({\rm det}X)^{p-|I|}{\rm
det}X_1(i_1)\cdots{\rm det}X_1(i_{|I|})=\frac{1}{|I|!}({\rm
det}X)^{p-1}{\rm det}X_1(I).$$ Dividing by $({\rm det}X)^{p-|I|}$
we obtain (ii).

By (ii) we have
\begin{align*}
\frac{1}{|I|!}({\rm
det}X)^{|I|-1}\frac{1}{|J|!}({\rm det}&X)^{|J|-1}{\rm
det}X_1(I){\rm det}X_1(J)=\\
&{\rm det}X_1(i_1)\cdots{\rm det}X_1(i_{|I|})
{\rm det}X_1(j_1)\cdots{\rm det}X_1(j_{|J|}).
\end{align*}
Thus if $\epsilon_{IJ}$ is the permutation arranging $I\cup J$ in
increasing order, then
$$
\frac{1}{|I|!|J|!} ({\rm
det}X)^{|I|+|J|-2} {\rm det}X_1(I) {\rm det}X_1(J) =
\frac{\epsilon_{IJ}}{(|I|+|J|)!} ({\rm det}X)^{|I|+|J|-1} {\rm
det}X_1(I\cup J).
$$
Now (iii) follows from dividing the above
equation by $({\rm
det}X)^{|I|+|J|-2}$ and multiplying by $|I|!|J|!$.
\end{proof}

Returning to the problem of finding the highest weight vector
associated to the Young diagram \eqnref{mainp>m}, our first task
is to present a more explicit expression for a product of the form
$\Gamma_{\la_{p+1}} \cdots \Gamma_{\la_{p+s}}$. We associate to
such a Young diagram a collection $D$ of $s$ marked diagrams
$D_i$, $i=1,\ldots,s$, with $D_i$ having $\la_{p+i}$ columns and
$m$ rows. We arrange these $D_i$ in the form:
\begin{equation}
{\beginpicture
\setcoordinatesystem units <1.5pc,1.5pc>
point at 0 5
\setplotarea x from -5 to 8, y from -4 to 12
\plot 0 8 6 8 6 12 0 12 0 8 /
\plot 0 9 6 9 /
\plot 0 10 6 10 /
\plot 0 11 6 11 /
\plot 1 8 1 12 /
\plot 2 8 2 12 /
\plot 3 8 3 12 /
\plot 4 8 4 12 /
\plot 5 8 5 12 /
\plot 0 3 5 3 5 7 0 7 0 3 /
\put{\vdots} at 1.5 1.5
\plot 0 -4 3 -4 3 0 0 0 0 -4 /
\put{$m\left\{\hbox to 0pt{
\vrule width 0 pc height 2.5pc depth 2.5pc\hfil}
  \right.$} at -0.6 10
\put{$m\left\{\hbox to 0pt{
\vrule width 0 pc height 2.5pc depth 2.5pc\hfil}
  \right.$} at -0.6 5
\put{$m\left\{\hbox to 0pt{
\vrule width 0 pc height 2.5pc depth 2.5pc\hfil}
  \right.$} at -0.6 -2
\put{$\underbrace{\hbox to 9pc{\ \hfill}}_{\la_{p+1}}$} at 3 7.5
\put{$\underbrace{\hbox to 7pc{\ \hfill}}_{\la_{p+2}}$} at 2.5 2.5
\put{$\underbrace{\hbox to 4pc{\ \hfill}}_{\la_{p+s}}$} at 1.5 -4.5
\put{$D_1:$} at -5 10
\put{$D_2:$} at -5 5
\put{$D_s:$} at -5 -2
\endpicture}
\end{equation}
Marked boxes are put into $D$ subject to the following constraint:
in each $D_i$ a column has at most one marked box. If a diagram
$D_i$ contains a marked box in its $k$-th row and $s$-th column,
then no other $D_j$ $(j\not =i)$ contains a marked box in its
$k$-th row and $s$-th column.  From now on $D$ will denote such a
collection of marked diagrams.

Now suppose that $D$ is a collection of diagrams $D_i$,
$i=1,\ldots,s$.  To each $D_i$ we may associate a
$\la_{p+i}\times\la_{p+i}$ matrix $Y_{D_i}$ as in the previous
section.  We let ${\rm det}Y_D:={\rm det}Y_{D_1}{\rm det}Y_{D_2}
\cdots {\rm det}Y_{D_s}$.  Now to each column $j$ of $D$ $(1\le
j\le \la_{p+1})$, we may associate a $p\times p$ matrix $X_j(I_j)$
obtained from the matrix $X$ as follows. Let $I_j$ be the subset
of $\{1, \ldots, p \}$ consisting of the numbers of the marked
rows on the column $j$. We define $X_j(I_j)$ to be the matrix
obtained from $X$ by replacing the rows of $X$ corresponding to
$I_j$ by the vector $(\eta^1_j, \ldots, \eta^p_j)$. We then define
$X_D:=X_1(I_1)X_2(I_2) \ldots X_{\la_{p+1}}(I_{\la_{p+1}})$.

Suppose we have a marked box in $D_i$ appearing in its $k$-th row
and $s$-th column.  We associate an odd indeterminate
$a_{i}^{ks}$.  Consider the product of all $a_i^{ks}$ arranged in
increasing order following the lexicographical ordering of
$(i,s,k)$.  Now we may also consider the product arranged in
increasing order following the lexicographical ordering $(s,k,i)$.
These two products differ by a sign, and this sign is denoted
by $\epsilon_D$.  Furthermore we let $d_i=|D_i|$, the number of
marked boxes in $D_i$, and $e_j$ be the number of marked
boxes in the $j$-th column of $D$.

\begin{prop}
With notation as above we have $$\frac{\Gamma_{\la_{p+1}} \cdots
\Gamma_{\la_{p+s}} }{({\rm det}X)^{\la_{p+2}+\ldots+\la_{p+s}}}=
\sum_{D}(-1)^{\frac{1}{2}((\sum_{i,j}d_id_j)-|D|)}
\frac{\epsilon_D}{e_1!\cdots e_{\la_{p+1}}!}{\rm det}X_D{\rm
det}Y_D.$$
\end{prop}
\begin{proof}  Given diagram $D_i$ with $m$ rows and $\la_{p+i}$ columns,
for $i=1,\ldots,s$, we want to know how to simplify the expression
$$(-1)^{\frac{1}{2}\sum_{i=1}^s d_i(d_i-1)}{\rm det}X_{D_1}{\rm
det}Y_{D_1}\cdots {\rm det}X_{D_s}{\rm det}Y_{D_s}.$$ We move all
${\rm det}X_{D_i}$ to the left and get
\begin{equation}\label{aux6}
(-1)^{\frac{1}{2}(\sum_{i,j=1}^s d_id_j)-\frac{1}{2}|D|}{\rm
det}X_{D_1}\cdots{\rm det}X_{D_s} {\rm det}Y_{D_1}\cdots{\rm
det}Y_{D_s}.
\end{equation}
Now each ${\rm det}X_{D_{i}}$ is a product of ${\rm det}X_b(a)$.
We arrange all $X_1(a)$ together so that $X_{1}(a_1)$ appears to
the left of $X_{1}(a_2)$ if and only if $a_1< a_2$ and call the
resulting expression $X_D^{\eta_1}$ and move it to the left.  We
do the same thing to $X_{2}(a)$ and move $X_D^{\eta_2}$ to the
right of $X_D^{\eta_1}$ etc. Then \eqnref{aux6} becomes
\begin{equation}\label{aux7}
(-1)^{\frac{1}{2}(\sum_{i,j=1}^s d_id_j)-\frac{1}{2}|D|}
\epsilon_DX_D^{\eta_1}X_D^{\eta_2}\cdots
X_D^{\eta_{\la_{p+1}}}{\rm det}Y_{D_1}\cdots{\rm det}Y_{D_s}.
\end{equation}
We apply now \corref{maincor} to \eqnref{aux7} and obtain
$$(-1)^{\frac{1}{2}(\sum_{i,j=1}^s d_id_j)-\frac{1}{2}|D|}
\frac{\epsilon_D}{e_1!\cdots e_{\la_{p+1}}!} ({\rm
det}X)^{(\la'_{p+1}-p)+\ldots+(\la'_{\la_{p+1}}-p)-\la_{p+1}} {\rm
det}X_D{\rm det}Y_D.$$ Since
$\sum_{i=1}^s\la_{p+i}=\sum_{i=1}^{\la_{p+1}}(\la'_{i}-p)$, the
proposition follows.
\end{proof}

It follows immediately from \corref{identity} that

\begin{prop} The vector
$$Z: ={\rm det}X_1(m+1,\ldots,p) \, {\rm det}
X_2(m+1,\ldots,p)\cdots {\rm det} X_{\la_{p+1}}(m+1,\ldots,p)$$ is
annihilated by \eqnref{radicalmn}.
\end{prop}

\begin{prop} The expression
$$ \sum_{D}(-1)^{\frac{1}{2}((\sum_{i,j}d_id_j)-|D|)}
\frac{\epsilon_D}{e_1!\cdots e_{\la_{p+1}}!} {\rm det}X_D {\rm
det}Y_D {\rm det}Z$$ is divisible by $({\rm det}X)^{\la_{p+1}}$.
Furthermore the resulting expression is independent of the
variables $x_l^i$ $(l =m+1,\ldots,p)$ and is annihilated by
\eqnref{radicalmn}.
\end{prop}
\begin{proof}
Since for a subset $I$ of $\{1,\ldots,m\}$, ${\rm det}X_i(I) {\rm
det}X_i(m+1,\ldots,p)$ is a scalar multiple of ${\rm det}X{\rm
det}X_i({I\cup\{m+1,\ldots,p\}})$ by \corref{maincor} (iii), it
follows that the expression is divisible by $({\rm
det}X)^{\la_{p+1}}$ and independent of $x_l^i$, for
$l=m+1,\ldots,p$, after division.  It is clear that it is
annihilated by \eqnref{radicalmn}.
\end{proof}

The vector $$ \sum_{D}(-1)^{\frac{1}{2}((\sum_{i,j}d_id_j)-|D|)}
(e_1!\cdots e_{\la_{p+1}}!)^{-1} \epsilon_D ({\rm det}X)^{-
\la_{p+1}}{\rm det}Z {\rm det}X_D {\rm det}Y_D $$ depends only on
the $s$-tuple $(\la_{p+1},\ldots,\la_{p+s})$ and thus we will
denote this vector by $\Gamma(\la_{p+1},\ldots,\la_{p+s})$.
\begin{prop}
The vector $\Gamma(\la_{p+1},\ldots,\la_{p+s})$ has weight
corresponding to the Young diagram \eqnref{mainp>m}.
\end{prop}

\begin{proof}
Let $f_j$, $j=1,\ldots,m$ denote the number of marked boxes in $D$
that appear in the $j$-th row of some diagram $D_i$.  Then ${\rm
det}Y_D$ has weight
$$\sum_{j=1}^{\la_{p+1}}(\la_j'-p)\delta_j+\sum_{j=1}^m
f_j\epsilon_j - \sum_{j=1}^{\la_{p+1}}e_j\delta_j
+\sum_{j=1}^{\la_{p+1}}(\la_j'-p)\tilde{\delta}_j,$$ while the
expression $({\rm det}X)^{- \la_{p+1}} {\rm det}X_D{\rm det}Z $
has weight
$$\la_{p+1}\sum_{j=1}^m\epsilon_j+(p-m)\sum_{j=1}^{\la_{p+1}}\delta_j
+\sum_{j=1}^{\la_{p+1}}e_j\delta_j-\sum_{j=1}^m
f_j\epsilon_j+\la_{p+1}\sum_{j=1}^p\tilde{\epsilon_j}.$$ So the
combined weight is
$$\la_{p+1}\sum_{j=1}^m\epsilon_j+\la_{p+1}\sum_{j=1}^p\tilde{\epsilon_j}
+\sum_{j=1}^{\la_{p+1}}(\la'_j-m)\delta_j+
\sum_{j=1}^{\la_{p+1}}(\la_j'-p)\tilde{\delta}_j,$$ which of
course is the weight of the Young diagram \eqnref{mainp>m}.
\end{proof}

Combining our results in this section we have proved

\begin{thm} In the case when $p\ge m$ an irreducible $gl(p|q)\times
gl(m|n)$ module $V_{p|q}^\la\otimes V_{m|n}^\la$ appears in
$\C[\x, \xibf, \etabf, \y]$ if and only if $\la_{m+1}\le n$ and
$\la_{p+1}\le q$.  The following are highest weight vectors
corresponding to such a Young diagram $\la$ ($t$ is the length of
$\la'$):
\begin{itemize}
\item[(i)] In the case when $\la_{m+1}=0$ it is given by
$$\prod_{i=1}^t\Delta_{\la'_i}.$$
\item[(ii)] In the case when $\la_{m+1}>0$ and $\la_{p+1}=0$ it is given by
$$\prod_{i=1}^r\Delta_{i,\la'_{i}}\prod_{i=r+1}^t\Delta_{\la'_{i}},$$
where $0\le r\le n$ is defined by $\la'_r>m$ and $\la'_{r+1}\le m$.
\item[(iii)] In the case when $\la_{p+1}>0$ it is given by
$$\Gamma(\la_{p+1},\ldots,\la_{p+s})\prod_{i=r+1}^{r'}\Delta_{i,\la'_{i}}
\prod_{i=r'+1}^t\Delta_{\la'_{i}},$$ where $r\le r'$ are defined
as in \eqnref{youngp>m} and $p+s$ is the length of $\la$.
\end{itemize}
\end{thm}

\begin{proof}
The only thing that remains to prove is that the vector in (iii)
is indeed killed by \eqnref{radicalmn}.  But this is because of
the presence of $Z$ in the formula of
$\Gamma(\la_{p+1},\ldots,\la_{p+s})$ and so is an immediate
consequence of \corref{identity}.
\end{proof}

\section{Construction of highest weight vectors in
$S({S^2\C^{m|n}})$}\label{hwvector}

In this section we will give an explicit formula for a highest
weight vector of each irreducible $gl(m|n)$-module that appear in
the symmetric algebra of the symmetric square of the natural
$gl(m|n)$-module.  According to \thmref{symmetric} we have the
following decomposition of $S({S^2\C^{m|n}})$ as a
$gl(m|n)$-module:
\begin{equation*}
S({S^2\C^{m|n}})\cong\sum_{\la}V_{m|n}^\la,
\end{equation*}
where the summation is over all partitions $\la$ with even rows
and $\la_{m+1}\le n$.

We let $\{x_1,\ldots,x_m;\xi_1,\ldots,\xi_n\}$ be the standard
basis for $\C^{m|n}$, with $x_i$ denoting even, and $\xi_j$ odd
vectors.  Regarding $x_i$ as even and $\xi_j$ as odd variables the
Lie supealgebra $gl(m|n)$ has a natural identification with the
space of first order differential operators over $\C[x_i,\xi_j]$.
The Cartan subalgebra of diagonal matrices is then spanned by its
standard basis $x_i\frac{\partial}{\partial x_i}$ and
$\xi_j\frac{\partial}{\partial\xi_j}$, for $i=1,\ldots,m$ and
$j=1,\ldots,n$.  The nilpotent radical is generated by the simple
root vectors
\begin{equation}\label{Borel}
x_i\frac{\partial}{\partial x_{i+1}},\quad
\xi_j\frac{\partial}{\partial \xi_{j+1}},\quad
x_m\frac{\partial}{\partial \xi_{1}},\qquad
i=1,\ldots,m-1;j=1,\ldots,n-1.
\end{equation}
${S^2\C^{m|n}}$ then is spanned by the vectors
$x_{ij}=x_{ji}=x_ix_j$, $y_{kl}=-y_{lk}=\xi_k\xi_l$ and
$\eta_{ki}=\xi_kx_i$, where $1\le i,j\le m$ and $1\le k,l\le n$.
This allows us to identify $S({S^2\C^{m|n}})$ with the polynomial
algebra over $\C$ in the even variables $x_{ij}$ and $y_{kl}$ and
odd variables $\eta_{ki}$, with $1\le i\le j\le m$ and $1\le
k<l\le n$, which we denote by $\C[x,y,\eta]$.

A convention of notation we will use throughout this section is
the following: By $x_i(x_1,x_2,\ldots,x_m)$ we will mean the row
vector $(x_{i1},x_{i2},\ldots,x_{im})$.  So by the expression
\begin{equation*}
\begin{pmatrix}
x_1(x_1,x_2,\ldots,x_m)\\ x_2(x_1,x_2,\ldots,x_m)\\ \vdots\\
x_m(x_1,x_2,\ldots,x_m)\\
\end{pmatrix}
\end{equation*}
we mean the matrix whose $i$-th row entries equals to
$(x_{i1},x_{i2},\ldots,x_{im})$, i.e.~the matrix
\begin{equation*}
X:=\begin{pmatrix} x_{11}&x_{12}&\cdots x_{1m}\\
x_{21}&x_{22}&\cdots x_{2m}\\ \vdots&\vdots&\cdots \vdots\\
x_{m1}&x_{m2}&\cdots x_{mm}\\
\end{pmatrix}.
\end{equation*}
Similarly by an expression of the form
\begin{equation*}
X_{i}(\xi_j):=\begin{pmatrix} x_1(x_1,x_2,\ldots,x_m)\\ \vdots\\
x_{i-1}(x_1,x_2,\ldots,x_m)\\
\xi_j(x_1,x_2,\ldots,x_m)\\x_{i+1}(x_1,x_2,\ldots,x_m)\\ \vdots\\
x_m(x_1,x_2,\ldots,x_m)\\
\end{pmatrix}
\end{equation*}
we mean to replace the $i$-th row of the matrix $X$ by the vector
$(\eta_{j1},\eta_{j2},\ldots,\eta_{jm})$. In these forms the
action of \eqnref{Borel} will be more transparent.

Consider the first $r\times r$ minor $\Delta_r$ of the $m\times m$
matrix $X$, for $1\le r\le m$.  It is easily seen to be a highest
weight vector in $\C[x,y,\eta]$ of highest weight
$2(\sum_{i=1}^r\epsilon_i)$, where as before we use $\epsilon_i$
and $\delta_k$ to denote the fundamental weights of $gl(m|n)$.
Hence if $\la$ is a Young diagram with even rows of length not
exceeding $m$, then its corresponding highest weight vector is a
product of $\Delta_r$'s.  To be explicit note that since $\la$ has
even rows, $\la_1=2t$ is an even number. Furthermore we also have
$\la'_{2i-1}=\la'_{2i}$, for all $i=1,\ldots,t$.  Then the highest
weight vector is given by $\prod_{i=1}^t\Delta_{\la'_{2i}}$.

Consider now a diagram of the form
\begin{equation}\label{basicdiagram}
{\beginpicture \setcoordinatesystem units <1.5pc,1.5pc> point at 0
3 \setplotarea x from 0 to 1.5, y from 0 to 5 \plot 0 0 4 0 /
\plot 0 0 0 5 / \plot 0 5 4 5 / \plot 4 5 4 0 / \plot 2 0 2 5 /
\plot 0 1 4 1 / \put{$m\left\{\hbox to 0pt{ \vrule width 0 pc
height 2.6pc depth 2.6pc\hfil}
  \right.$} at -.5 3
\put{$1\left\{\hbox to 0pt{ \vrule width 0 pc height 0.5pc depth
0.5pc\hfil}
  \right.$} at -.4 .5
\put{$\underbrace{\hbox to 6pc{\ \hfill}}_{2l}$} at 2 -0.5
\endpicture}
\end{equation}
The product of highest weight vectors of two such diagrams, if
non-zero, gives a highest weight vector of a diagram of the form
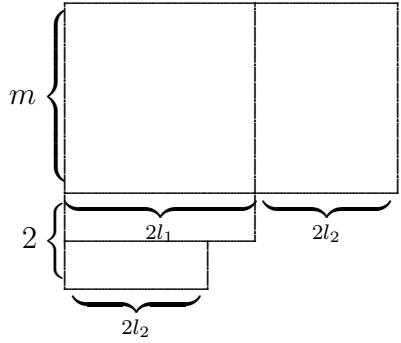
\begin{equation*}
{\beginpicture \setcoordinatesystem units <1.5pc,1.5pc> point at 0
3 \setplotarea x from 0 to 1.5, y from 0 to 6 \plot 0 1 4 1 /
\plot 0 1 0 6 / \plot 0 6 4 6 / \plot 4 6 4 1 / \plot 0 2 4 2 /
\plot 4 2 7 2 / \plot 7 2 7 6 / \plot 7 6 4 6 / \plot 0 0 0 1 /
\plot 0 0 3 0 / \plot 3 0 3 1 / \put{$m\left\{\hbox to 0pt{ \vrule
width 0 pc height 2.6pc depth 2.6pc\hfil}
  \right.$} at -.5 4
\put{$2\left\{\hbox to 0pt{ \vrule width 0 pc height 1pc depth
1pc\hfil}
  \right.$} at -.4 1
\put{$\underbrace{\hbox to 6pc{\ \hfill}}_{2l_1}$} at 2 1.4
\put{$\underbrace{\hbox to 4pc{\ \hfill}}_{2l_2}$} at 5.5 1.4
\put{$\underbrace{\hbox to 4pc{\ \hfill}}_{2l_2}$} at 1.5 -.6
\endpicture}
\end{equation*}
Dividing by the determinant ${\rm det}X^{l_2}$ we obtain a highest
weight vector for the diagram
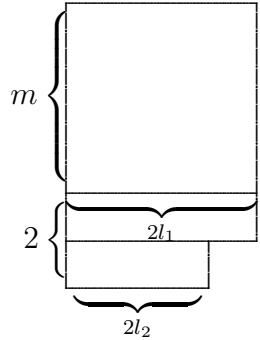
\begin{equation*}
{\beginpicture \setcoordinatesystem units <1.5pc,1.5pc> point at 0
3 \setplotarea x from 0 to 1.5, y from 0 to 6 \plot 0 1 4 1 /
\plot 0 1 0 6 / \plot 0 6 4 6 / \plot 4 6 4 1 / \plot 0 2 4 2 /
\plot 0 0 0 1 / \plot 0 0 3 0 / \plot 3 0 3 1 /
\put{$m\left\{\hbox to 0pt{ \vrule width 0 pc height 2.6pc depth
2.6pc\hfil}
  \right.$} at -.5 4
\put{$2\left\{\hbox to 0pt{ \vrule width 0 pc height 1pc depth
1pc\hfil}
  \right.$} at -.4 1
\put{$\underbrace{\hbox to 6pc{\ \hfill}}_{2l_1}$} at 2 1.4
\put{$\underbrace{\hbox to 4pc{\ \hfill}}_{2l_2}$} at 1.5 -.6
\endpicture}
\end{equation*}
Thus it is enough to construct vectors associated to the Young
diagrams of the form \eqnref{basicdiagram}.  To do so we first
consider the case when $l=1$ in \eqnref{basicdiagram}.

Consider the expression
\begin{align}\label{delta12}
\Delta(\xi_1,\xi_2):=&-({\rm det}X)(\xi_1\xi_2)+({\rm
det}X_1(\xi_1))(\xi_2x_1)\\ &+({\rm
det}X_2(\xi_1))(\xi_2x_2)+\cdots+({\rm
det}X_m(\xi_1))(\xi_2x_m),\nonumber
\end{align}
where by $(\xi_1\xi_2)$ and $(\xi_2x_i)$ we mean $y_{12}$ and
$\eta_{2i}$, respectively.  The following lemma will be useful
later on.

\begin{lem}\label{auxilary}
Let $A=(a_{ij})$ be a complex symmetric $m\times m$ matrix and
$\theta_1,\theta_2,\ldots,\theta_m$ be odd variables.  Then
\begin{equation*}
{\rm det}\begin{pmatrix} 0&\theta_1\cdots\theta_m\\ \theta_1&\\
\vdots&A\\ \theta_m&\\
\end{pmatrix}=0.
\end{equation*}
\end{lem}

\begin{proof}
It is enough to restrict ourselves to real symmetric $m\times m$
matrices $A$.  Let $U$ be an orthogonal $m\times m$ matrix such
that $U^{t}AU=D$, where $D$ is a diagonal matrix.  We compute
\begin{equation}\label{aux33}
\begin{pmatrix}
1&0\cdots 0\\ 0&\\ \vdots&U^t\\ 0&\\
\end{pmatrix}
\begin{pmatrix}
0&\theta_1\cdots\theta_m\\ \theta_1&\\ \vdots&A\\ \theta_m&\\
\end{pmatrix}
\begin{pmatrix}
1&0\cdots 0\\ 0&\\ \vdots&U\\ 0&\\
\end{pmatrix}=
\begin{pmatrix}
0&\zeta_1\cdots\zeta_m\\ \zeta_1&\\ \vdots&D\\ \zeta_m&\\
\end{pmatrix},
\end{equation}
where $\zeta_k=\sum_{j=1}^mu_{jk}\theta_j$ and $U=(u_{ij})$.  But
the determinant of the matrix on the right-hand side of
\eqnref{aux33} is zero.
\end{proof}

The next lemma is straightforward.

\begin{lem}
$\Delta(\xi_1,\xi_2)$ has weight
$2(\sum_{i=1}^m\epsilon_i)+\delta_1+\delta_2$ and hence its weight
corresponds to the weight of the Young diagram
\eqnref{basicdiagram} with $l=1$.
\end{lem}

\begin{lem}\label{lem1}
$\Delta(\xi_1,\xi_2)$ is annihilated by all operators in
\eqnref{Borel} and hence is a highest weight vector in
$S({S^2\C^{m|n}})$ corresponding to the Young diagram
\eqnref{basicdiagram} with $l=1$.
\end{lem}

\begin{proof}
First consider the action of the operator
$x_{i-1}\frac{\partial}{\partial x_{i}}$, for $i=2,\ldots,m$, on
$\Delta(\xi_1,\xi_2)$ given as in \eqnref{delta12}.  Certainly
$x_{i-1}\frac{\partial}{\partial x_{i}}$ annihilates the first
summand of \eqnref{delta12}, and furthermore it takes the summand
$X_j(\xi_1)(\xi_2x_j)$ for $j\not=i-1,i$ to
\begin{equation*}
{\rm det}\begin{pmatrix}
x_1(x_1,\ldots,x_{i-1},x_{i-1},\ldots,x_m)\\ \vdots\\
x_{j-1}(x_1,\ldots,x_{i-1},x_{i-1},\ldots,x_m)\\
\xi_1(x_1,\ldots,x_{i-1},x_{i-1},\ldots,x_m)\\
x_{j+1}(x_1,\ldots,x_{i-1},x_{i-1},\ldots,x_m)\\ \vdots\\
x_m(x_1,\ldots,x_{i-1},x_{i-1},\ldots,x_m)
\end{pmatrix}(\xi_2x_j)+
{\rm det}\begin{pmatrix} x_1(x_1,\ldots,x_m)\\ \vdots\\
x_{j-1}(x_1,\ldots,x_m)\\ \xi_1(x_1,\ldots,x_m)\\
x_{j+1}(x_1,\ldots,x_m)\\ \vdots\\ x_{i-1}(x_1,\ldots,x_m)\\
x_{i-1}(x_1,\ldots,x_m)\\ \vdots\\ x_m(x_1,\ldots,x_m)
\end{pmatrix}(\xi_2x_j),
\end{equation*}
which is zero. $x_{i-1}\frac{\partial}{\partial x_{i}}$ takes
$X_i(\xi_1)(\xi_2x_i)$ to
\begin{equation*}
{\rm det}\begin{pmatrix}
x_1(x_1,\ldots,x_{i-1},x_{i-1},\ldots,x_m)\\ \vdots\\
x_{i-i}(x_1,\ldots,x_{i-1},x_{i-1},\ldots,x_m)\\
\xi_1(x_1,\ldots,x_{i-1},x_{i-1},\ldots,x_m)\\
x_{i+1}(x_1,\ldots,x_{i-1},x_{i-1},\ldots,x_m)\\ \vdots\\
x_m(x_1,\ldots,x_{i-1},x_{i-1},\ldots,x_m)
\end{pmatrix}(\xi_2x_i)+
{\rm det}\begin{pmatrix} x_1(x_1,\ldots,x_m)\\ \vdots\\
x_{i-1}(x_1,\ldots,x_m)\\ \xi_1(x_1,\ldots,x_m)\\
x_{i+1}(x_1,\ldots,x_m)\\ \vdots\\ x_m(x_1,\ldots,x_m)
\end{pmatrix}(\xi_2x_{i-1}).
\end{equation*}
The first summand is zero, while the second summand remains. Now
we verify similarly that $x_{i-1}\frac{\partial}{\partial x_{i}}$
takes $X_{i-1}(\xi_1)(\xi_2x_{i-1})$ to the identical expression
as the second summand above with the difference that the $i-1$-st
and $i$-th rows are interchanged.  Thus
$x_{i-1}\frac{\partial}{\partial x_{i}}(\Delta(\xi_1,\xi_2))=0$.

Consider now the action of $x_m\frac{\partial}{\partial\xi_1}$ on
$\Delta(\xi_1,\xi_2)$.  Note that
$x_m\frac{\partial}{\partial\xi_1}$ kills every term in
\eqnref{delta12} except for the first and the last.  The
contribution from the first summand is $-{\rm det}X(\xi_2x_m)$,
while that from the last summand is ${\rm det}X(\xi_2x_m)$, and
hence $x_m\frac{\partial}{\partial\xi_1} (\Delta(\xi_1,\xi_2))=0$.

Finally we consider the action of
$\xi_1\frac{\partial}{\partial\xi_2}$ on $(\Delta(\xi_1,\xi_2))$
as in \eqnref{delta12}.  $\xi_1\frac{\partial}{\partial\xi_2}$
kills the first term in \eqnref{delta12} and the resulting vector
is
\begin{equation}\label{aux22}
\sum_{i=1}^m ({\rm det}X_{i}(\xi_1))(\xi_1x_i),
\end{equation}
which can be in a consistent form with our earlier notation
written as $\Delta(\xi_1,\xi_1)$. Expanding along the first row we
see that \eqnref{aux22} is the same as
\begin{equation*}
{\rm det}\begin{pmatrix} 0&(\xi_1x_1)\cdots(\xi_1x_m)\\
(\xi_1x_1)&\\ \vdots&X\\ (\xi_1x_m)&\\
\end{pmatrix},
\end{equation*}
which is zero by \lemref{auxilary}.
\end{proof}

The proof of the above theorem gives us certain identities that
will be used later on.  We will collect them here for the
convenience of the reader:
\begin{align}
x_{i-1}\frac{\partial}{\partial
x_i}(\Delta(\xi_k,\xi_l))&=0,\label{id1}\\
\xi_{j-1}\frac{\partial}{\partial
\xi_j}(\Delta(\xi_j,\xi_l))&=\Delta(\xi_{j-1},\xi_l),\label{id12}\\
\xi_{j-1}\frac{\partial}{\partial
\xi_j}(\Delta(\xi_s,\xi_j))&=\Delta(\xi_s,\xi_{j-1}),\label{id3}\\
\Delta(\xi_j,\xi_j)&=0.\label{id4}
\end{align}

We now turn our attention to the general case of a Young diagram
of the form \eqnref{basicdiagram} with general $l$.  Of course we
have the restriction that $2l\le n$.

Let $\sigma=\{(i_1,i_2),(i_3,i_4),\ldots,(i_{2l-1},i_{2l})\}$ be a
partition of the set $\{1,2,\ldots,2l\}$.  Assuming that we have
arranged $\sigma$ in the form so that
$i_1<i_2,i_3<i_4,\ldots,i_{2l-1}<i_{2l}$, we may define
$\epsilon_{\sigma}$ to be the sign of the permutation taking $k$
to $i_k$ for all $1\le k\le 2l$. We may associate to $\sigma$ a
vector
$\Delta(\xi_{i_1},\xi_{i_2})\cdots\Delta(\xi_{i_{2l-1}},\xi_{i_{2l}})$
in $S({S^2\C^{m|n}})$ and define
\begin{equation*}
\Gamma(2l)=\sum_{\sigma}\epsilon_\sigma\Delta(\xi_{i_1},\xi_{i_2})
\cdots\Delta(\xi_{i_{2l-1}},\xi_{i_{2l}}),
\end{equation*}
where the sum is taken over all partition
$\sigma=\{(i_1,i_2),(i_3,i_4),\ldots,(i_{2l-1},i_{2l})\}$ of the
set $\{1,2,\ldots,2l\}$ arranged in the form that
$i_1<i_2,i_3<i_4,\ldots,i_{2l-1}<i_{2l}$.  The following lemma is
again a straightforward computation.

\begin{lem}
The weight of $\Gamma(2l)$ is
$2l(\sum_{i=1}^m\epsilon_i)+\sum_{j=1}^{2l}\delta_j$ and hence
corresponds to the weight of the Young diagram
\eqnref{basicdiagram}.
\end{lem}

\begin{lem}
$\Gamma(2l)$ is annihilated by \eqnref{Borel} and hence a highest
weight vector in $S({S^2\C^{m|n}})$.
\end{lem}

\begin{proof}
The fact that $\Gamma(2l)$ is annihilated by
$x_{i-1}\frac{\partial}{\partial x_i}$ for $i=2,\ldots,m$ is a
consequence of \eqnref{id1}.  Now the proof of \lemref{lem1} shows
that $x_m\frac{\partial}{\partial\xi_1}(\Delta(\xi_1,\xi_s))=0$,
for every $s=1,\ldots,n$.  Thus
$x_m\frac{\partial}{\partial\xi_1}$ annihilates $\Gamma(2l)$ as
well. So it remains to show that
$\xi_{j-1}\frac{\partial}{\partial\xi_j}$ kills $\Gamma(2l)$.

Given a summand in $\Gamma(2l)$ of the form
$\epsilon_{\sigma}\cdots\Delta(\xi_{j-1},\xi_{l})\cdots\Delta(\xi_j,\xi_i)\cdots$
there exists a summand of the form
$\epsilon_{\sigma'}\cdots\Delta(\xi_{j-1},\xi_{i})\cdots\Delta(\xi_j,\xi_l)\cdots$,
which is identical to it except at these two places.  Now
$\xi_{j-1}\frac{\partial}{\partial\xi_j}$ takes the first of the
two summands above to
$$\epsilon_{\sigma}\cdots\Delta(\xi_{j-1},\xi_{l})\cdots\Delta(\xi_{j-1},\xi_i)\cdots,$$
and the second summand to
$$\epsilon_{\sigma'}\cdots\Delta(\xi_{j-1},\xi_{i})\cdots\Delta(\xi_{j-1},\xi_l)\cdots.$$
But $\sigma$ and $\sigma'$ differ by a transposition $(i,l)$ and
hence $\epsilon_\sigma=-\epsilon_{\sigma'}$ and so these two terms
cancel.

Consider a summand in $\Gamma(2l)$ of the form
$\epsilon_{\sigma}\cdots\Delta(\xi_{l},\xi_{j-1})\cdots\Delta(\xi_j,\xi_i)\cdots$.
But in $\Gamma(2l)$ we also have a summand of the form
$\epsilon_{\sigma'}\cdots\Delta(\xi_{l},\xi_{j})\cdots\Delta(\xi_{j-1},\xi_i)\cdots$.
Applying $\xi_{j-1}\frac{\partial}{\partial\xi_j}$ to these two
terms, we again see that they cancel by the same reasoning as
before.

Now we look at a term of the form
$\epsilon_{\sigma}\cdots\Delta(\xi_{j},\xi_{l})\cdots\Delta(\xi_i,\xi_{j-1})\cdots$.
We also have a term of the form
$\epsilon_{\sigma'}\cdots\Delta(\xi_{j-1},\xi_{l})\cdots\Delta(\xi_{i},\xi_j)\cdots$.
Again they will cancel each other after applying
$\xi_{j-1}\frac{\partial}{\partial\xi_j}$.

Finally a term of the form
$\epsilon_{\sigma'}\cdots\Delta(\xi_{j-1},\xi_{j})\cdots$ is
killed by $\xi_{j-1}\frac{\partial}{\partial\xi_j}$ by
\eqnref{id4}.  This completes the proof.
\end{proof}

It is clear that a product of $\Gamma(2l)$'s (not necessarily the
same $l$) is non-zero, which therefore allows us to construct all
other highest weight vectors, as discussed in the beginning of
this section.  Below we
summarize the results of this section.

\begin{thm}\label{hwv}
The $gl(m|n)$-highest weight vectors of $S({S^2\C^{m|n}})$ form an
abelian semigroup generated by
$\Gamma(2),\ldots,\Gamma(2[\frac{n}{2}])$ and
$\Delta_1,\ldots,\Delta_m$, where $[\frac{n}{2}]$ denotes the
largest integer not exceeding $\frac{n}{2}$.  Furthermore this
semigroup is free if and only if $n=0,1$.  More precisely a
highest weight vector associated to an even partition
$\la=(\la_1,\ldots,\la_l)$ with $\la_{m+1}\le n$ is given by
\begin{equation*}
{({\rm det}X)^{-\sum_{i=m+2}^l\la_i}}\prod_{m+1}^l
\Gamma(\la_i)\prod_{j=r+1}^l(\Delta_{\la'_j})^{\frac{1}{2}},
\end{equation*}
where the non-negative integer $r$ is defined by $\la'_{r}>m$ and
$\la'_{r+1}\le m$.
\end{thm}

\begin{rem}\label{generic}
From \thmref{hwv} we may recover the highest weight vectors in
$S(S^2\C^m)$ and $S(\Lambda^2\C^n)$ by putting $n=0$ and $m=0$,
respectively.  Namely, identifying $S^2\C^m$ (respectively
$\Lambda^2\C^n$) with the space of symmetric $m\times m$
(respectively skew-symmetric $n\times n$) matrices, we see that in
the first case the highest weight vectors are generated by the
leading minors of the determinant of the typical element of
$S^2\C^m$, while in the second case they are generated by the
Pfaffians of the leading $2l\times 2l$ minors of the the typical
element of $\Lambda^2\C^n$, where $2l\le n$. (cf.~\cite{H2}).
\end{rem}

\end{document}